# A theoretical framework and some promising findings of grey wolf optimizer, part II: global convergence analysis

Haoxin Wang, Libao Shi*
National Key Laboratory of Power Systems, Shenzhen International Graduate School, Tsinghua University, 518055, P. R. China

## Abstract

This paper proposes a theoretical framework of the grey wolf optimizer (GWO) based on several interesting theoretical findings, involving sampling distribution, order-1 and order-2 stability, and global convergence analysis. In the part II of the paper, the global convergence analysis is carried out based on the well-known stagnation assumption for simplification purposes. Firstly, the global convergence property of the GWO under stagnation assumption is abstracted and modelled into two propositions, corresponding to global searching ability analysis and probability-1 global convergence analysis. Then, the global searching ability analysis is carried out. Next, based on a characteristic of the central moments of the new solution of the GWO under stagnation assumption, the probability-1 global convergence property of the GWO under stagnation assumption is proved. Finally, all conclusions are verified by numerical simulations, and the corresponding discussions that the global convergence property can still be guaranteed in the original GWO without stagnation assumption are given.

## 1. Introduction

In recent decades, a kind of optimization algorithms with random searching mechanism, called meta-heuristic algorithms, have been elaborately designed and developed to solve some complicated optimization problems existed in real life [1-2]. As one of the typical meta-heuristic algorithms, the Grey Wolf Optimizer (GWO) [3] has attracted widespread attention, due to its intuitive program and excellent optimization performance. So far, most of the research conducted on the GWO has mainly focused on algorithmic variants and engineering applications [4-14], and the theoretical analysis on algorithm performance requires in-depth research.

In the part I of the paper, as fundamental tools for performance analysis of the GWO, the sampling distribution characteristics of the new solution were elaborately discussed, after which the order-1 and order-2 stability under stagnation assumption were proved. In the part II of the paper, based on the findings in the part I, the global convergence properties of the GWO under stagnation assumption will be discussed in depth from theoretical analysis and numerical simulations.

It is known that the global convergence analysis is an important part of random search meta-heuristic algorithm performance analysis. Early work in convergence analysis mainly focused on several convergence criteria called algorithm condition and convergence condition proposed by Solis F. [15], which will be elaborated in section 5. Based on these convergence criteria, the (1+1) evolution strategy (ES) with only a constant scale factor mutation operator was proved to be a global convergence algorithm, while the global convergence of (1+1) ES with a self-adaptive scale factor mutation operator was no more guaranteed [16-19]. In addition, it was discovered that the global convergence of particle swarm optimizer (PSO) was not guaranteed [20-21], and a benchmark function example on which the differential evolution (DE) failed to find the global optimum was provided and tested by theoretical analysis and numerical simulations [22-23]. Later work mainly focused on how to offer new tools to model the convergence problem and how to analyze the convergence properties of some novel meta-heuristic algorithms or variants of classical meta-heuristic algorithms. For example, the global convergence of the PSO was studied further by Markov chain [24-25], martingale theory [26], probabilistic metric space [27], and spectral analysis [28], etc. The global convergence properties of two novel meta-heuristics called modified fireworks explosion optimization (mFEO) [29] and clouds search optimization (CSO) [30] were studied in [31], and the convergence performance of several variants of the PSO (such as quantum PSO [32-33], MOPSO [34-35], and bare bones PSO [36-37], etc.) were discussed as well.

In the part II of the paper, the global convergence properties of the GWO are proved and discussed elaborately. Similar to the part I, the stagnation assumption is taken into account throughout the analysis of this paper for simplification. Firstly, two propositions are proposed for the global searching ability analysis and the probability-1 global convergence analysis respectively, and then the two propositions are proved and carefully discussed, during which several properties of the sampling distribution of the new solution in the GWO under stagnation assumption as discussed in the part I are reintroduced as crucial lemmas. Finally, numerical simulations are carried out to testify the conclusions drawn in the part II of the paper. It should be noticed that the global searching ability of the GWO is substantially reduced under the stagnation assumption, thus the global convergence of the GWO without stagnation assumption can be assumed.

The part II of the paper is organized as follows. In section 2, the summarized work pertinent to part I is introduced briefly. In section 3, the global convergence properties of the GWO under stagnation assumption are abstracted into two propositions corresponding to the global searching ability analysis and the probability-1 global convergence analysis, and the global searching ability analysis is conducted. In section 4, the probability-1 global convergence property of the GWO under stagnation assumption is proved. In section 5, the corresponding numerical simulations and discussions are provided. In section 6, we conclude the results above and point out the limitations of our research.

## 2. A summary of part I

To facilitate the discussions in the following sections, several main conclusions drawn in the part I of the paper are summarized in this section. The solution updating mechanism of the GWO is shown as:

$$\mathbf{x_i}(t+1) = \frac{1}{3}\sum_{k=1}^{3}\mathbf{x'_k}(t) \tag{2.1}$$

$$x'_{kj}(t) = p_{kj}(t) + A_{kj}\left|C_{kj}p_{kj}(t) - x_{ij}(t)\right| \tag{2.2}$$

$$\mathbf{p_1}(t+1) = \begin{cases}\mathbf{x_i}(t+1), f\big(\mathbf{x_i}(t+1)\big) < f\big(\mathbf{p_1}(t)\big)\\ \mathbf{p_1}(t), \text{otherwise}\end{cases} \tag{2.3}$$

$$\mathbf{p_2}(t+1)=\begin{cases}\mathbf{x_i}(t+1), f(\mathbf{p_1})\le f(\mathbf{x_i}(t+1))<f(\mathbf{p_2}(t))\\ \mathbf{p_2}(t), \text{otherwise}\end{cases} \tag{2.4}$$

$$\mathbf{p_3}(t+1)=\begin{cases}\mathbf{x_i}(t+1), f(\mathbf{p_2})\le f(\mathbf{x_i}(t+1))<f(\mathbf{p_3}(t))\\ \mathbf{p_3}(t), \text{otherwise}\end{cases} \tag{2.5}$$

where all variables mentioned in this section are explained in Table I, and the procedure of the GWO can be found in Fig. 1 of the part I.

Table I Nomenclature

| Variable type | Variable | Meaning |
|---|---|---|
| Problem information | $\min f(\mathbf{x})$ | Optimization problem to be solved |
| | $D$ | Dimension of the problem |
| Algorithm parameters | $T$ | Total iterations |
| | $N$ | Number of search agents |
| Indices | $t$ | Index of iteration |
| | $i$ | Index of search agent |
| | $j$ | Index of dimension |
| Solutions | $\mathbf{p_1}(t)=[p_{1j}(t)\mid j=1,\dots,D]$, $\mathbf{p_2}(t),\mathbf{p_3}(t)$ | Positions of the best three search agents at iteration $t$ |
| | $\mathbf{x_i}(t)=[x_{ij}(t)]$ | Position of the $i$th search agent at iteration $t$ |
| Parameters in updating solution | $\mathbf{x'_k}(t)=[x'_{kj}(t)], k=1,2,3$ | Intermediate vector at iteration $t$ |
| | $A_{kj}$ | Random variable following $U[-a,a]$ |
| | $C_{kj}$ | Random variable following $U[0,2]$ |
| | $a$ | $2\left(1-\frac{t}{T}\right)$ |
| Intermediate variables proposed for simplification | $m_{kj}$ | $a(-\lvert p_{kj}(t)\rvert+\lvert x_{ij}(t)-p_{kj}(t)\rvert)$ |
| | $n_{kj}$ | $a(\lvert p_{kj}(t)\rvert+\lvert x_{ij}(t)-p_{kj}(t)\rvert)$ |
| | $\sum p_j, \sum m_j, \sum n_j$ | $\sum_{k=1}^{3} p_{kj}(t), \sum_{k=1}^{3} m_{kj}, \sum_{k=1}^{3} n_{kj}$ |
| PDFs | $g_{kj}(u_j)$ | PDF of $x'_{kj}(t)$ |
| | $h_{ij}(u_j)/h_{ij,t+1}(u_j)$ | PDF of $x_{ij}(t+1)$ |

In the part I of the paper, the sampling distribution that $\mathbf{x_i}(t)$ follows was elaborately analyzed. Under the stagnation assumption that all $\mathbf{p_k}(t)(k=1,2,3)$ were constant vectors independent of $t$ and the constant $\mathbf{x_i}(t)$ assumption, the PDF of $x_{ij}(t+1)$ (marked as $h_{ij}(u_j)$) can be expressed as the convolution of the PDFs of $x'_{kj}(t)$ ($k=1,2,3$, marked as $g_{kj}(u_j)$). Further discussions showed that the domain of $h_{ij}(u_j)$ was $\left(\frac{1}{3}\sum p_j-\frac{1}{3}\sum n_j, \frac{1}{3}\sum p_j+\frac{1}{3}\sum n_j\right)$, and $h_{ij}(u_j)$ was a single-peak function symmetrical about $u_j=\frac{1}{3}\sum p_j$. After eliminating the constant $\mathbf{x_i}(t)$ assumption, the domain of the PDF of $x_{ij}(t)$ (marked as $h_{ij,t}(u_j)$) became an open interval containing $\left(\frac{1}{3}\sum p_j-\max_{x_{ij}(1)}\left|x_{ij}(1)-\frac{1}{3}\sum p_j\right|\prod_{\tau=1}^{t-1}2\left(1-\frac{\tau}{T}\right), \frac{1}{3}\sum p_j+\max_{x_{ij}(1)}\left|x_{ij}(1)-\frac{1}{3}\sum p_j\right|\prod_{\tau=1}^{t-1}2\left(1-\frac{\tau}{T}\right)\right)$ for $t\ge 2$, however, $h_{ij,t}(u_j)$ was still a single-peak function symmetrical

about $u_j = \frac{1}{3}\sum p_j$. In addition, the recursive expression of the central moments of any positive integer order of $x_{ij}(t)$ was derived elaborately, and finally, the order-1 and order-2 stability of the GWO under stagnation assumption were proved.

# 3. Mathematical modelling of GWO global convergence analysis under stagnation assumption

In this section, the global convergence of the GWO is modelled and abstracted into two propositions, corresponding to the global searching ability analysis and the probability-1 global convergence analysis. It should be noticed that the *stagnation assumption* mentioned in section 2 will be taken into account throughout the whole discussions of the global convergence analysis of the GWO in following sections. However, in section 5.2, the impact of the stagnation assumption on the global convergence property of the original GWO will be discussed. Similar to the part I, since the positions of the best three search agents $\mathbf{p_k}(t)$ remain unchanged throughout iterations, in this section and section 4, the $\mathbf{p_k}(t)$ will be abbreviated as $\mathbf{p_k}$.

## 3.1 Problem formulation

In the part II of this paper, the optimization problem to be solved is marked as $\min f(\mathbf{x})$. We assume that the global optimum of $f(\mathbf{x})$ is marked as $\mathbf{x_g}$, and a set $S_\varepsilon$ is defined by the fitness value as:

$$S_\varepsilon = \{\mathbf{x} | f(\mathbf{x}) - f(\mathbf{x_g}) \le \varepsilon\} \tag{3.1}$$

In addition, the $f(\mathbf{x})$ is assumed to be a continuous function with only one global optimal solution $\mathbf{x_g}$ (i.e. $\forall \mathbf{x} \in \mathbb{R}^D \backslash \{\mathbf{x_g}\}$ , $f(\mathbf{x_g}) < f(\mathbf{x})$), in which case the $S_\varepsilon$ defined by Eq. (3.1) becomes a neighborhood of the global optimal solution $\mathbf{x_g}$ as long as the value of $\varepsilon$ is small enough. Fig. 1 (a) illustrates the case that the value of $\varepsilon$ is small enough, and Fig. 1 (b) illustrates the case that the value of $\varepsilon$ is large.

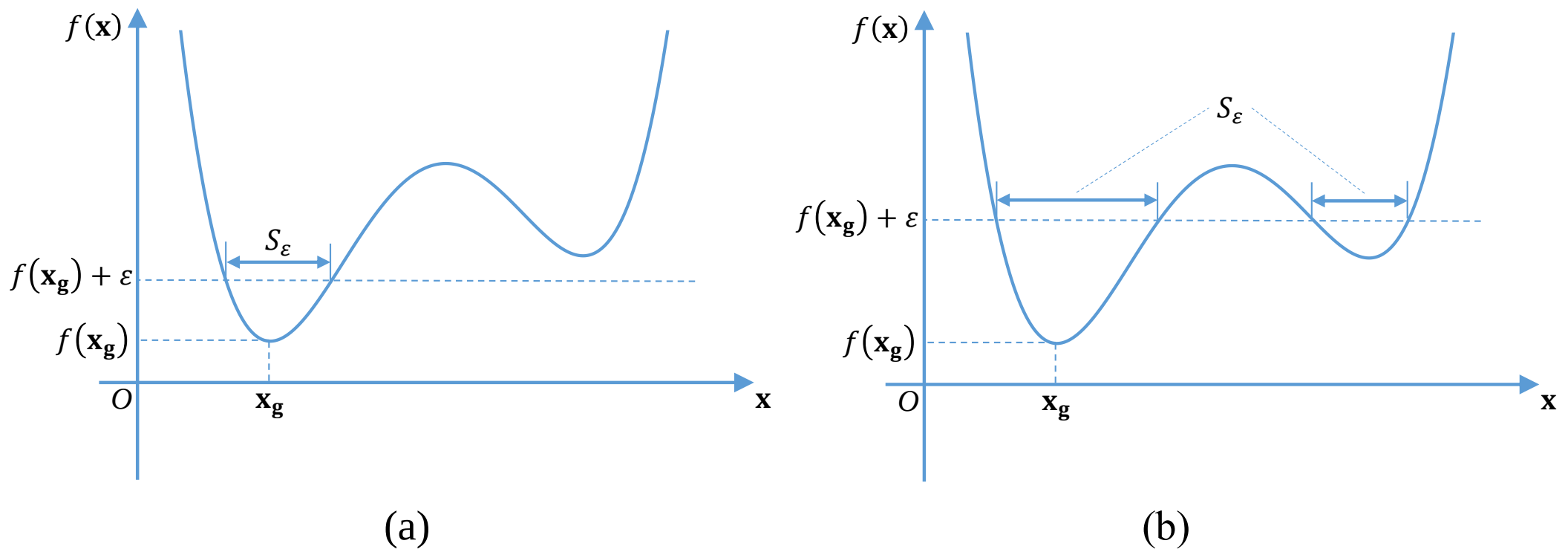

Fig. 1 $S_\varepsilon$ pertinent to a multi-modal one-dimension function $f(\mathbf{x})$, (a): small value $\varepsilon$, (b): large value $\varepsilon$.

Since $f(\mathbf{x})$ is a continuous function, there must exist a box region

$$S_0 = \{\mathbf{x} = [x_j] | lb_j \le x_j \le ub_j, j = 1, \dots, D, lb_j < ub_j\} \quad (3.2)$$

with positive Lebesgue measure satisfying $S_0 \subseteq S_\varepsilon$. Since the Lebesgue measure of the exact global optimum set $\{\mathbf{x_g}\}$ is zero, which means the exact position of $\mathbf{x_g}$ will never be found via any random search algorithm, thus in this paper, if any solution within $S_0$ (a positive Lebesgue measure set) is found, we will say that the global optimum is found. For instance, regarding the optimization problem $\min_{\mathbf{x}\in\mathbb{R}^2} f(\mathbf{x}) = x_1^2 + 2x_2^2$ with minimum value 0 and $\mathbf{x_g} = [0,0]$, it can be solved that $S_\varepsilon = \{(x_1, x_2) | x_1^2 + 2x_2^2 \le \varepsilon\}$ is within an ellipse, and we can choose $S_0 = \left\{(x_1, x_2) | -\frac{\sqrt{\varepsilon}}{\sqrt{2}} \le x_1 \le \frac{\sqrt{\varepsilon}}{\sqrt{2}}, -\frac{\sqrt{\varepsilon}}{2} \le x_2 \le \frac{\sqrt{\varepsilon}}{2}\right\}$.

## 3.2 Mathematical modelling of global convergence analysis of the GWO

In this section, the mathematical model of the global convergence analysis of the GWO under stagnation assumption is established.

If an optimization algorithm holds global convergence, it can be easily understood that for each search agent the ability to search the global optimum (i.e. $S_0$ in Eq. (3.2)) is always guaranteed throughout iterations. As for the GWO under stagnation assumption, there exist two barriers in successfully finding the global optimum:

(1) In the GWO under stagnation assumption, the initial positions of all search agents are randomly chosen within the feasible region, thus the global optimum $\mathbf{x_g}$ may not be within the envelope of the initial positions of all search agents (for example, regarding the following case: the feasible region is $[-1,1]^D$, the global optimum $\mathbf{x_g} = [0.8, \dots, 0.8]$, and the initial positions of $N$ agents are randomly chosen within the feasible region, then the probability that all of the agents are within $[-1,0.6]^D$ becomes $\left(\frac{0.6-(-1)}{1-(-1)}\right)^{ND} = 0.8^{ND} > 0$, in which case $\mathbf{x_g}$ must goes beyond the envelope of all initial positions). According to the part I, the region of $x_{ij}(t+1)$ is bounded within $\left(\frac{1}{3}\sum p_j - \frac{1}{3}\sum n_j, \frac{1}{3}\sum p_j + \frac{1}{3}\sum n_j\right)$, thus, is it true that for each search agent the ability to find the global optimum is guaranteed after a few iterations $t_0$ (i.e. $\mathbb{P}\{x_{ij}(t_0+1) \in S_0\} > 0$), or the ability to find the global optimum is no longer guaranteed when the initial positions of search agents are specially chosen (i.e. $\forall t, \mathbb{P}\{x_{ij}(t) \in S_0\} = 0$)?

(2) According to the part I, the mutation step size of each solution in the GWO under stagnation assumption is proportional to $a(t)$, which decreases during iterations. This means that the mutation step size of each solution tends to decrease during iterations, as shown in Fig. 9 of the part I. Thus, after enough iterations, is each search agent able to find the global optimum with probability 1? Or when the mutation step size decreases sharply, does a search agent fail to find the optimum (i.e. find a position within $S_0$ defined in Eq. (3.2), similarly hereinafter) with probability 1?

Fortunately, it can be proved that these two barriers are able to overcome. The corresponding conclusions can be abstracted into the following two propositions, respectively:

**Proposition 1**

In the GWO under stagnation assumption, $\forall \mathbf{p_k} \in \mathbb{R}^D, \mathbf{x_i}(1) \neq \frac{1}{3}\sum_{k=1}^{3} \mathbf{p_k}, \forall S_0 \subseteq \mathbb{R}^D$ , $\exists T, t_0 \leq \left\lfloor \frac{T}{2} \right\rfloor \in N^*, s.t. \forall t = t_0, t_0 + 1, \ldots, \left\lfloor \frac{T}{2} \right\rfloor, \mathbb{P}\{\mathbf{x_i}(t) \in S_0\} > 0.$

**Proposition 2**

In the GWO under stagnation assumption, $\forall \mathbf{p_k} \in \mathbb{R}^D, \mathbf{x_i}(1) \neq \frac{1}{3}\sum_{k=1}^{3} \mathbf{p_k}, \lim_{t\to\infty} \mathbb{P}\{\mathbf{x_i}(t) \in S_0\} = 1.$

Obviously, the proposition 1 points out that the global searching ability of each agent is always guaranteed no matter how the initial position is chosen, thus the analysis pointing towards proposition 1 is called *global searching ability analysis*. The proposition 2 points out that the probability of each agent to successfully find the global optimum tends to 1 after enough iterations, thus the analysis pointing towards proposition 2 is called *probability-1 global convergence analysis*. Since the probability that the initial positions of all search agents are identical with each other (i.e. $\mathbf{x_i}(1) = \frac{1}{3}\sum_{k=1}^{3} \mathbf{p_k}$) is 0, thus by proving the two propositions above, it can be asserted that the GWO under stagnation assumption holds global convergence. The proof of proposition 1 can be found in Appendix (the proofs of all lemmas given below can be found in Appendix), and the proposition 2 will be discussed and proved in section 4.

Finally, the proposition 1 is exhibited in Fig. 2, in which the green point denotes the mutation center $\frac{1}{3}\sum_{k=1}^{3} \mathbf{p_k}$ (unchanged during iterations under the stagnation assumption), the purple point denotes the initial position of search agent $i$ ($\mathbf{x_i}(1)$), the big nesting cuboids $D_i(t)$ denote the searching region of $\mathbf{x_i}(t)$, and the cuboid containing the orange point denotes the neighborhood $S_0$ of the global optimum $\mathbf{x_g}$. It has been proved in the part I (please see the proposition 3.1 in the part I) that the region of $x_{ij}(t)$ contains $\left(\frac{1}{3}\sum p_j - \max_{x_{ij}(1)} \left| x_{ij}(1) - \frac{1}{3}\sum p_j \right| \prod_{\tau=1}^{t-1} 2\left(1 - \frac{\tau}{T}\right), \frac{1}{3}\sum p_j + \max_{x_{ij}(1)} \left| x_{ij}(1) - \frac{1}{3}\sum p_j \right| \prod_{\tau=1}^{t-1} 2\left(1 - \frac{\tau}{T}\right)\right)$, thus in the early stage of the algorithm when $a > 1$, the searching region of $\mathbf{x_i}(t)$ expands during iterations, thus able to contain any point within the feasible region (such as $\mathbf{x_g}$) after a few iterations.

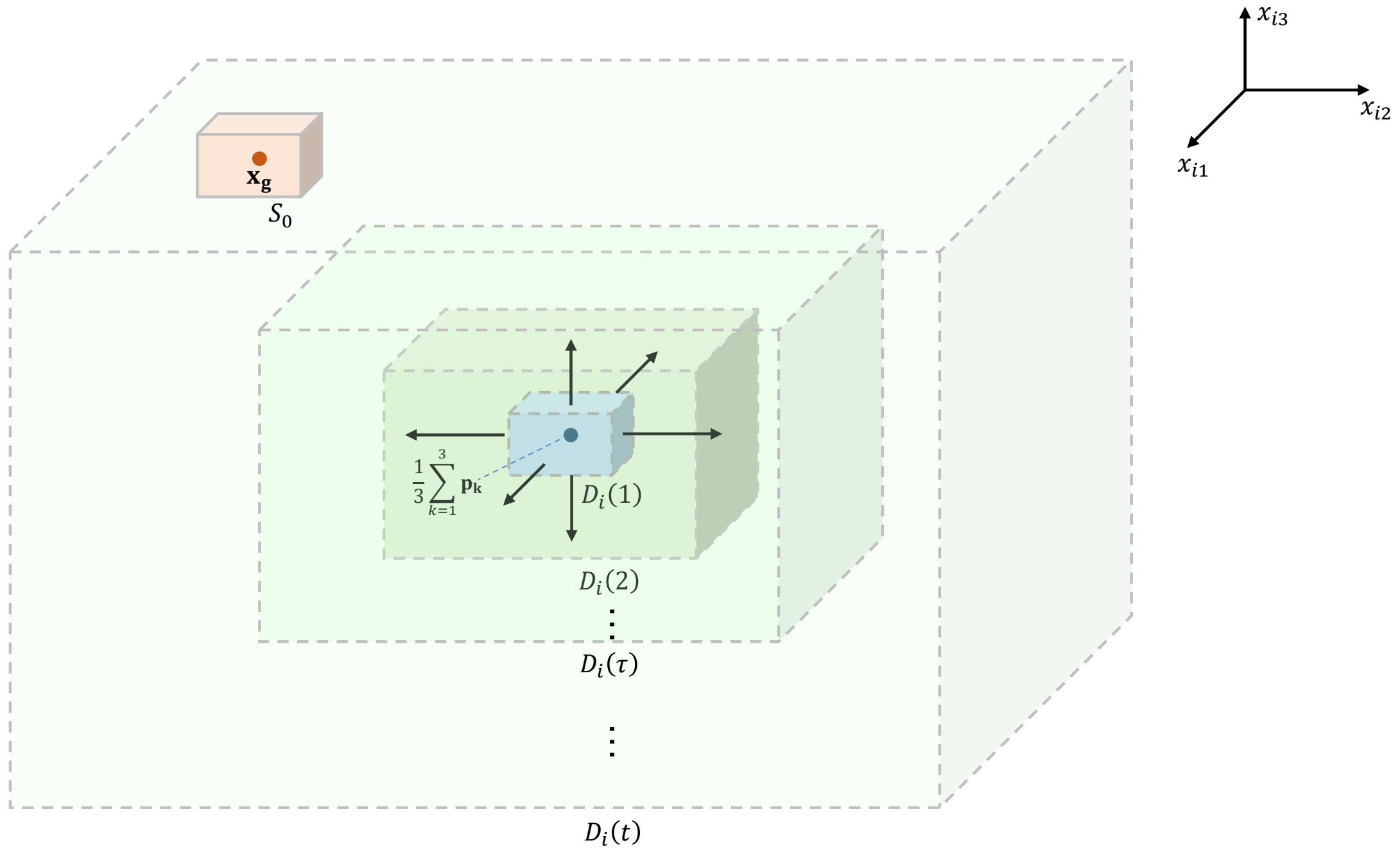

Fig. 2 Global searching process of $\mathbf{x_i}(t)$

# 4. Probability-1 global convergence analysis of the GWO under stagnation assumption

In this section, the probability-1 global convergence analysis of the GWO under stagnation assumption is carried out, and the proposition 2 proposed in section 3.2 is proved and carefully discussed.

## 4.1 Introduction of algorithm condition and convergence condition

Since $\mathbf{x_i}(t)$ is a complicated stochastic process, it's almost impossible to calculate the value of $\mathbb{P}\{\mathbf{x_i}(t) \in S_0\}$ directly. Fortunately, as mentioned in section 1, the studies on several convergence criteria [15] of random searching optimization algorithms can serve as fundamental analysis tools. In [15], in order to solve the optimization problem $\min f(\mathbf{x})$, an optimization algorithm based on random searching mechanism is defined as:

$$\mathbf{z}(t+1) = D\big(\mathbf{z}(t), \boldsymbol{\xi}(t)\big) \tag{4.1}$$

where $\mathbf{z}(t)$ denotes the solution at iteration *t*, $\boldsymbol{\xi}(t)$ denotes the new sample generated by probability measure $\mu_t$, and *D* denotes the operator to generate the new solution $\mathbf{z}(t+1)$ from the original solution $\mathbf{z}(t)$ and new sample $\boldsymbol{\xi}(t)$.

In addition, it is pointed out in [15] that:

**Lemma 1**

For the solution of an unconstrained optimization problem, a random searching optimization algorithm defined by Eq. (4.1) converges to the global optimum with probability 1 if the following two conditions are satisfied:

(1) (Algorithm condition):

$$f(\mathbf{z}(t+1)) \le \min\{f(\mathbf{z}(t)), f(\boldsymbol{\xi}(t))\} \tag{4.2}$$

(2) (Convergence condition): For any set $A$ with positive Lebesgue measure,

$$\prod_{t=1}^{\infty}(1-\mu_t(A)) = 0 \tag{4.3}$$

where $\mu_t(A)$ is the probability of $A$ being generated by $\mu_t$.

The specific meanings of the algorithm condition and convergence condition are introduced as follows: In the algorithm condition, $f(\mathbf{z}(t+1)) \le f(\mathbf{z}(t))$ is required, which means the fitness curve of the algorithm should be a non-increasing function of iteration $t$. In addition, $f(\mathbf{z}(t+1)) \le f(\boldsymbol{\xi}(t))$ is required, which means if an individual $\boldsymbol{\xi}(t)$ with better fitness value than $\mathbf{z}(t)$ is found, then at least the position of $\boldsymbol{\xi}(t)$ should be adopted as $\mathbf{z}(t+1)$, which means each effective search (i.e. a solution with better fitness value is found) will not be wasted. As for the convergence condition, Eq. (4.3) denotes that the value of $\mu_t(A)$ should be large enough (if $\forall t, \mu_t(A) = 0$, then $\prod_{t=1}^{\infty}(1-\mu_t(A)) = 1$ holds), i.e., for any set $A$ with positive measurement, the probability that a new solution is found within $A$ cannot be too small, which indicates that the global searching ability of the algorithm is always guaranteed. Please refer to [15] for more details about these two conditions.

In the GWO under stagnation assumption, the $\mathbf{z}(t)$ in Eq. (4.1) means the best search agent $\mathbf{p_1}(t)$, $\boldsymbol{\xi}(t)$ means the new solution $\mathbf{x_i}(t+1)$, and $\mu_t$ means the probability distribution that $\mathbf{x_i}(t+1)$ follows. According to the stagnation assumption, the position of the best search agent remains unchanged until the global optimum ($S_0$) is found, thus the updating mechanism of $\mathbf{p_1}(t+1)$ $D$ becomes:

$$\mathbf{p_1}(t+1) = \begin{cases} \mathbf{x_i}(t+1), \mathbf{x_i}(t+1) \in S_0 \\ \mathbf{p_1}(t), \text{otherwise} \end{cases} \tag{4.4}$$

Therefore, in order to prove proposition 2, we only need to prove that the GWO under stagnation assumption satisfies both algorithm condition and convergence condition. Since $\mathbf{p_1}(t) \notin S_0$ (otherwise the algorithm succeeds in finding the global optimum in iteration 1), thus when $\mathbf{x_i}(t+1) \in S_0, f(\mathbf{x_i}(t+1)) \le f(\mathbf{p_1}(t))$ holds, accordingly we have:

**Lemma 2**

The GWO under stagnation assumption satisfies the algorithm condition in lemma 1.

Thus, we only need to prove that the GWO under stagnation assumption satisfies the convergence condition, which will be elaborated in section 4.3.

## 4.2 Central moments of $x_{ij}(t)$ revisited

It has been concluded in the part I that the distribution of $x_{ij}(t)$ becomes more close to $\frac{1}{3}\sum p_j$ during iterations, which indicates the decrease of the global searching ability of agent $i$. However, if the value of the maximum iteration time $T$ is large enough, then in the early stage of iteration (for example, from $t = 0$ to $t = 0.05T$), the speed that the value of $a(t)$ decreases is small enough that the global searching ability of agent *i* is always guaranteed during this stage, in which case we can optimistically suppose that when $T \to \infty$, the probability that agent *i* finds the global optimum may tend to 1, just as the proposition 2 points out. Evidently, in order to analyze the global searching ability of search agent *i* quantitatively, more information of the distribution that $x_{ij}(t)$ follows is required, in which the results of the moments of $x_{ij}(t)$ as derived in the part I of the paper provide the corresponding theoretical analysis basis. Therefore, in this section, we will further discuss the characteristics of the moments of $x_{ij}(t)$. As a core conclusion of this section, lemma 4 will serve as a crucial lemma in the proof of proposition 2.

It has been derived in the part I that when *r* is an even number, the *r*-order central moment of $x_{ij}(t)$ (marked as $\sigma^r x_{ij}(t)$) can be expressed in the recurrence form below:

$$\sigma^r x_{ij}(t+1) = f_r(t)\sigma^r x_{ij}(t) + g_r(t), r \text{ is an even number} \tag{4.5}$$

where

$$f_r(t) = \frac{a^r}{3^r}\sum_{p=0}^{\frac{r}{2}}\sum_{q=0}^{\frac{r}{2}-p}\frac{\binom{r}{2p}\binom{r-2p}{2q}}{(2p+1)(2q+1)(r-2p-2q+1)}, r \text{ is an even number} \tag{4.6}$$

and $g_r(t)$ denotes a linear combination of $\sigma^m x_{ij}(t)(m = 0,1,\dots,r-1)$. In addition, when *r* is an odd number, we have

$$\sigma^r x_{ij}(t) = 0, r \text{ is an odd number} \tag{4.7}$$

The following conclusion pertinent to $f_r(t)$ is crucial in subsequent discussions:

**Lemma 3**

(1) $f_r(t) = \frac{a^r}{3^r}\frac{3^{r+3}-3}{4(r+1)(r+2)(r+3)}$

(2) Solve the root of equation $f_r(t) = 1$ about $a$ as (marked as $a_r$, similarly hereinafter):

$$a_r = 3\left(\frac{4(r+1)(r+2)(r+3)}{3^{r+3}-3}\right)^{\frac{1}{r}} \tag{4.8}$$

then the value of $a_r$ decreases with the increase of positive even number *r*.

(3) when $r$ is an even number and $a \in [1.9,2]$, $f_r(t) < 1$ holds when $r \le 6$, and $f_r(t) > 1$ holds when $r \ge 8$.

When $r = 2,4,6,8,10$, the values of $a_r$ calculated from Eq. (4.8) are given in Table II respectively, and the solutions of $f_r(t) \ge 1$ are shown in Fig. 3, where the green arrow represents the iteration process of the GWO under stagnation assumption (i.e. $a$ decreases from 2 to 0).

Table II Values of $a_r$ when $r = 2,4,6,8,10$, respectively

| $r$ | 2 | 4 | 6 | 8 | 10 |
|---|---|---|---|---|---|
| $a_r$ | 3 | 2.36 | 2.05 | 1.87 | 1.74 |

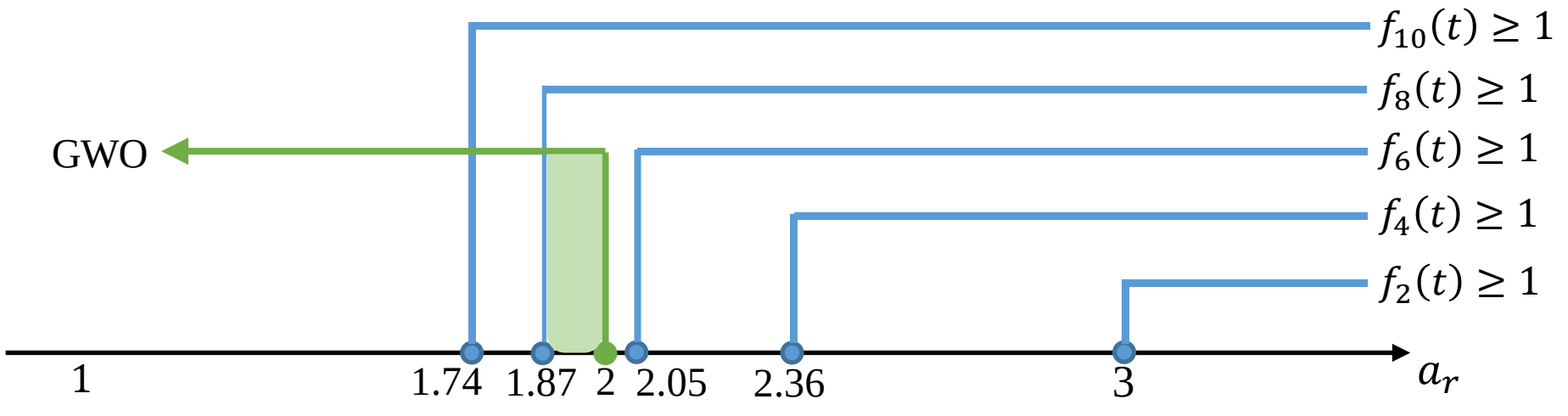

Fig. 3 Solutions of $f_r(t) \geq 1$ when $r = 2{,}4{,}6{,}8{,}10$, respectively

It can be seen from Fig. 3 that $a_r$ is a decreasing function of positive even number *r*, and $\{a|f_2(t) \geq 1\} \subset \{a|f_4(t) \geq 1\} \subset \{a|f_6(t) \geq 1\} \subset \cdots$ holds. In addition, since $a_8 \approx 1.87 < 1.9 < 2 < a_6 \approx 2.05$, the conclusion in lemma 3(3) is exhibited directly.

In the following part of section 4, we will limit the discussion to the early stage of the GWO under stagnation assumption, which is defined as $t \leq \lfloor 0.05T \rfloor$. In this case, the value of $a$ is bounded within $[1.9{,}2]$, and $a_8 \approx 1.87 < a < a_6 \approx 2.05$ holds. According to the algorithm condition, the optimum found before iteration *t* (i.e. the fitness value curve) will not increase with *t*, thus the probability-1 global convergence property of the GWO under stagnation will become evident if the probability-1 global convergence property can be proved during $t \leq \lfloor 0.05T \rfloor$. Please notice that when $\sigma^m x_{ij}(t)(m = 0{,}1, \dots, r-1)$ are solved, the values of $f_r(t)$ and $g_r(t)$ are determined, and Eq. (4.5) denotes a first-order system of $\sigma^r x_{ij}(t)$.

We suppose that when $a \in [1.9{,}2]$, $\forall m = 0{,}1, \dots, r-1, \lim_{t \to \infty} \sigma^m x_{ij}(t) < \infty$ exists (since $\sigma^1 x_{ij}(t) = 0$, with the assumption that $\sigma^0 x_{ij}(t) = 1$, then this assumption holds when $r = 2$, which means this assumption is reasonable), thus $\lim_{t \to \infty} g_r(t) < \infty$ exists. Evidently, the system defined by Eq. (4.5) converges to a certain value when $0 < f_r(t) < 1$ holds for any $a \in [1.9{,}2]$, and diverges when $f_r(t) > 1$ holds for any $a \in [1.9{,}2]$. In addition, since $\forall a \in [1.9{,}2], a_2 > a_4 > a_6 > a > a_8 > a_{10} > \cdots$, and $f_2(t), f_4(t), f_6(t) < 1, f_8(t), f_{10}(t), \dots > 1$, holds, we can draw the following conclusion, i.e. the lemma4 which will serve as a crucial lemma in the proof of proposition 2 in section 4.3:

**Lemma 4**

In the GWO under stagnation assumption, for any positive even number $r \leq 6, \lim_{T \to \infty} \lim_{t \to \lfloor 0.05T \rfloor} \sigma^r x_{ij}(t) < +\infty$ exists, while for any positive even number $r \geq 8$, $\lim_{T \to \infty} \lim_{t \to \lfloor 0.05T \rfloor} \sigma^r x_{ij}(t) = +\infty$.

**Remark 1**

It has been proved in the part I of the paper that for any fixed $a \in [0{,}2]$, the variance of $x_{ij}(t)$ ($\sigma^2 x_{ij}(t)$) converges when $t \to \infty$. This conclusion can be testified according to the discussions mentioned above as well. According to Eq. (4.6), $f_2(t) = \frac{a^2}{9}\sum_{p=0}^{1}\sum_{q=0}^{1}\frac{\binom{2}{2p}\binom{2-2p}{2q}}{(2p+1)(2q+1)(3-2p-2q)} = \frac{a^2}{9} < 1$ holds, since $\sigma^0 x_{ij}(t) = 1$ and $\sigma^1 x_{ij}(t) = 0$ then $\lim_{t\to\infty}\sigma^2 x_{ij}(t)$ exists.

## 4.3 Convergence condition for the GWO under stagnation assumption

In this section, we will prove that the GWO under stagnation assumption satisfies the convergence condition described in lemma 1.

According to the convergence condition as given in Eq. (4.3), the item $\left(1-\mu_t(A)\right)$ denotes the probability that the algorithm fails to find a solution within $A$. Taking the arbitrariness of $\mathbf{p_k}, \mathbf{x_i}(1)$ and $S_0 = \{\mathbf{x}|lb_j \le x_j \le ub_j\}$ into account, the set $A$ can be substituted by $S_0$. In addition, since there are $N$ search agents in the swarm, the $\left(1-\mu_t(A)\right)$ can be reformulated as:

$$\left(1-\mu_t(A)\right) = \left(1-\mu_t(S_0)\right) = \mathbb{P}\{\mathbf{x_i}(t+1) \notin S_0 | \forall i = 1,\dots,N\} = \prod_{i=1}^{N}\mathbb{P}\{\mathbf{x_i}(t+1) \notin S_0\}$$

$$= \prod_{i=1}^{N}(1-\mathbb{P}\{\mathbf{x_i}(t+1) \in S_0\})$$

where the second equation is obtained according to the definition of the GWO under stagnation assumption, and the third equation can be obtained because all $\mathbf{x_i}(t+1)(i=1,\dots,N)$ are independent of each other. In addition, since all $x_{ij}(t+1)(j=1,\dots,D)$ are independent of each other, and $S_0 = \{\mathbf{x}|lb_j \le x_j \le ub_j\}$ is a box set, the item $\left(1-\mu_t(A)\right)$ can be further reformulated as:

$$\left(1-\mu_t(A)\right) = \prod_{i=1}^{N}(1-\mathbb{P}\{\mathbf{x_i}(t+1) \in S_0\}) = \prod_{i=1}^{N}\left(1-\prod_{j=1}^{D}\mathbb{P}\{x_{ij}(t+1) \in [lb_j, ub_j]\}\right)$$

Consequently, in order to prove that the GWO under stagnation assumption satisfies the convergence condition, we only need to prove that:

**Lemma 5**

$$\prod_{t=1}^{\infty}\prod_{i=1}^{N}\left(1-\prod_{j=1}^{D}\mathbb{P}\{x_{ij}(t+1) \in [lb_j, ub_j]\}\right) = 0 \quad (4.9)$$

It has been mentioned in the part I of the paper that the PDF of $x_{ij}(t+1)$ (marked as $h_{ij,t}(u_j)$) cannot be derived directly, thus the item $\mathbb{P}\{x_{ij}(t+1) \in [lb_j, ub_j]\}$ cannot be derived analytically. Fortunately, the item $\mathbb{P}\{x_{ij}(t+1) \in [lb_j, ub_j]\}$ can be roughly approximated based on the information of the central moments of $x_{ij}(t+1)$ as calculated and analyzed in section 4.2.

Let us reconsider the discussions in section 4.2. According to lemma 4, when $T \to \infty$ and $t \to \lfloor 0.05T \rfloor$, it can be concluded that the $x_{ij}(t)$ will approach to the probabilistic distribution of a random variable (marked as $x_{ij}(\infty)$) whose value of central moment $\sigma^r x_{ij}(\infty)$ is finite when $r \le 6$ or infinite when $r \ge 8$. Although the expression of $h_{ij,t}(u_j)$ cannot be derived precisely, the characteristics of the limit distribution of $x_{ij}(t)$ (whose PDF is marked as $h_{ij,\infty}(u_j)$) can be analyzed by calculating the value of $\sigma^r x_{ij}(\infty)$ when $T \to \infty$ and $t \to \lfloor 0.05T \rfloor$. In this section, the lemma 6 is designed to indicate the following characteristic of $h_{ij,\infty}(u_j)$:

**Lemma 6**

$\forall u_0 \in \mathbb{R}, h_{ij,\infty}(u_0) = 0$ is a false proposition.

The meaning of lemma 6 is obvious: in order to draw the conclusion of global convergence, the global searching ability of each agent should always be guaranteed, i.e., the speed that the searching range of each agent decreases should be slow enough. According to lemma 6, for $\forall u_0 \in \mathbb{R}$, the value of $h_{ij,t}(u_0)$ will not tend to 0 when $T \to \infty$ and $t \to \lfloor 0.05T \rfloor$, which means for $\forall u_0$, the searching ability of $x_{ij}(t+1)$ around $u_0$ is guaranteed, i.e. the corresponding global searching ability is guaranteed.

Lemma 6 provides a characteristic of the PDF $h_{ij,\infty}(u_0)$, which can be easily rewritten in the following probability form:

**Lemma 7**

$\lim_{T\to\infty} \lim_{t\to\lfloor 0.05T \rfloor} \mathbb{P}\{x_{ij}(t+1) \in [lb_j, ub_j]\} = 0$ is a false proposition.

Lemma 7 provides a more intuitive conclusion of the global searching ability of $x_{ij}(t+1)$, that is to say, no matter how the initial position of each search agent is chosen, the global searching ability of $x_{ij}(t+1)$ will not tend to 0 when $T \to \infty$ and $t \to \lfloor 0.05T \rfloor$, in another word, the global searching ability is always guaranteed throughout iterations. Since lemma 5 is a simple corollary of lemma 7 (please see the Appendix), the proof of the global convergence of the GWO under stagnation assumption with probability 1 is finished.

In summary, the procedure in proving the global convergence of the GWO under stagnation assumption with probability 1 is shown in Fig. 4

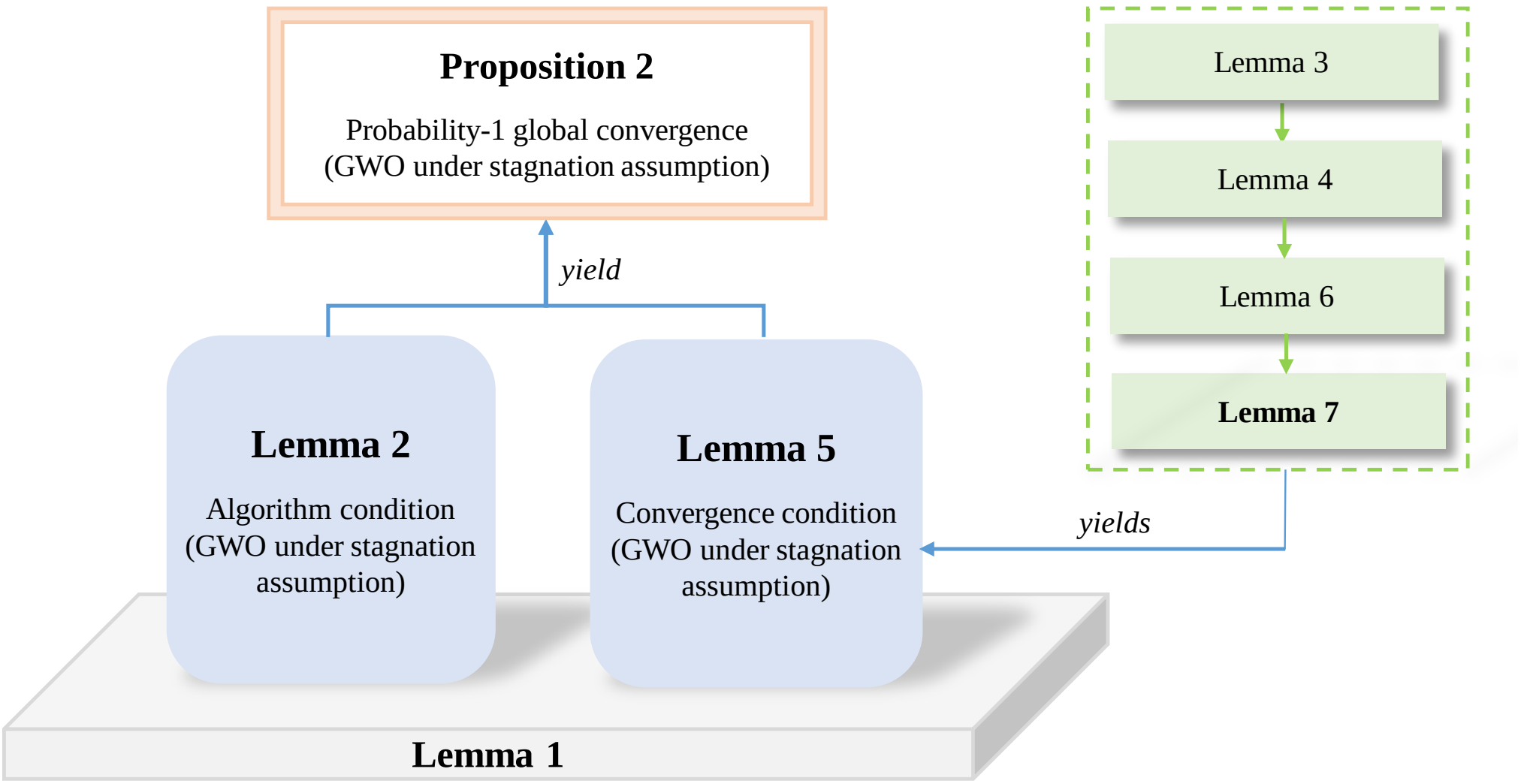

Fig. 4 Procedure in proving the probability-1 global convergence of GWO under stagnation assumption

# 5. Numerical simulations and discussions

In this section, numerical simulations are carried out to verify the theoretical conclusions drawn in section 4. In addition, the impact of stagnation assumption on the global convergence property of the original GWO will be discussed. It can be seen from both quantitative analysis and numerical simulations that the stagnation assumption reduces the convergence property of the original GWO sharply, thus the global convergence of the original GWO can be asserted as well.

## 5.1 Numerical simulations

In this section, the probability-1 global convergence properties of the GWO will verified by numerical simulations. However, instead of directly verifying the proposition 2, the lemmas given in section 4 will be verified in this section, which is because the probability-1 global convergence of the GWO under stagnation assumption cannot be verified directly by applying the GWO under stagnation assumption to solve various test problems. Specifically speaking, if the GWO under stagnation assumption fails to solve an optimization problem, we cannot determine the failure of finding the global optimum is due to whether the global convergence of the GWO under stagnation assumption is not guaranteed or the total number of iterations and the number of search agents are not large enough. If the GWO under stagnation assumption is able to converge to the global optimum (as proved above), then no matter what the test problem is, the GWO under stagnation assumption will be able to find the global optimum with probability 1 as long as the total number of iterations and the number of search agents are large enough, which are impractical for simulations since the computational cost is limited.

**Remark 2**

In order to better understand the difficulty in verifying the probability-1 global convergence of the GWO under stagnation assumption by applying the GWO under stagnation assumption to solve various test problems, the following optimization problem $\min y = \sum_{j=1}^{D}\left(x_j - 1\right)^2, -10 \le x_j \le 10$ with global optimum $\mathbf{x_g} = [1, \dots, 1]$ is solved as benchmark. The number of search agents is set to be $N = 10$, the maximum number of iterations $T$ is set to be a sufficiently large number (such that the GWO under stagnation assumption will succeed in finding the global optimum before iteration $T$), the GWO under stagnation assumption is terminated when a solution is found within a sphere neighborhood of $\mathbf{x_g}$: $S_\varepsilon = \{\mathbf{x} | f(\mathbf{x}) - f(\mathbf{x_g}) < 0.1\}$, and the positions of the best three search agents are set as $\mathbf{p_1} = [3, \dots, 3], \mathbf{p_2} = [4, \dots, 4], \mathbf{p_3} = [5, \dots, 5]$ (these positions remain unchanged during iterations under the stagnation assumption). For different values of $D$, the average minimum number of iterations before successfully finding the global optimum ($S_\varepsilon$) is recorded in Table II:

Table II Relationship between computational cost and $D$

| $D$ | 2 | 3 | 4 | 5 |
|---|---|---|---|---|
| average minimum number of | 104 | 5992 | 374610 | 29457387 |
| iterations | | | | |
| average minimum calculation time ($s$) | 0.009 | 0.42 | 25.75 | 2116.04 |

According to Table II, when $D$ increases, the average minimum number of iterations increases sharply. It should be noted that the neighborhood of $\mathbf{x_g}$ defined by $S_\varepsilon = \{\mathbf{x} | f(\mathbf{x}) - f(\mathbf{x_g}) < 0.1\}$ is quite loose, in another word, it is hardly practicable to verify the probability-1 global convergence of the GWO under stagnation assumption by applying the GWO under stagnation assumption to solve various test problems.

Since the lemma 1 has been carefully discussed in [15], the lemmas 2-3 are obvious and easy to test, and the lemma 5 is a direct conclusion of lemma 7, thus, in this section, only the lemma 4, lemma 6, and lemma 7 will be verified by numerical simulations successively.

### 5.1.1 Verification of the lemma 4

In order to verify the lemma 4, the GWO under stagnation assumption is performed from $a = 2$ to $a = 1.9$, and the values of $\sigma^r x_{ij}(t)(r = 2,4,6,8,10,12)$ are simulated, respectively. Here, we set $T = 2000$ (since the value of $T$ is large enough, the condition $T \to \infty$ is assumed to be satisfied, similarly hereinafter), and the values of $p_{kj}$ are randomly chosen. For each group of $p_{kj}$ chosen, the GWO under stagnation assumption is run for $10^6$ trials from $t = 1$ to $t = \lfloor 0.05T \rfloor = 100$. If the value of $x_{ij}$ obtained in $l$th trial at iteration $t$ is marked as $x_{ij}(l,t)$, and $\forall l, x_{ij}(l,1)$ is chosen as a random variable following $U[-4,4]$, then the value of $\sigma^r x_{ij}(t)$ can be approximated as:

$$\sigma^r x_{ij}(t) \approx \frac{1}{10^6 - 1} \sum_{l=1}^{10^6} \left( x_{ij}(l,t) - 10^{-6} \sum_{q=1}^{10^6} x_{ij}(q,t) \right)^r \quad (5.1)$$

In addition, considering that the range of $\sigma^r x_{ij}(t)$ is sometimes relatively large, the approximate value of $\ln \sigma^r x_{ij}(t)$ is shown instead of $\sigma^r x_{ij}(t)$. The corresponding simulation results are shown in Fig. 5 (a)-(f), respectively.

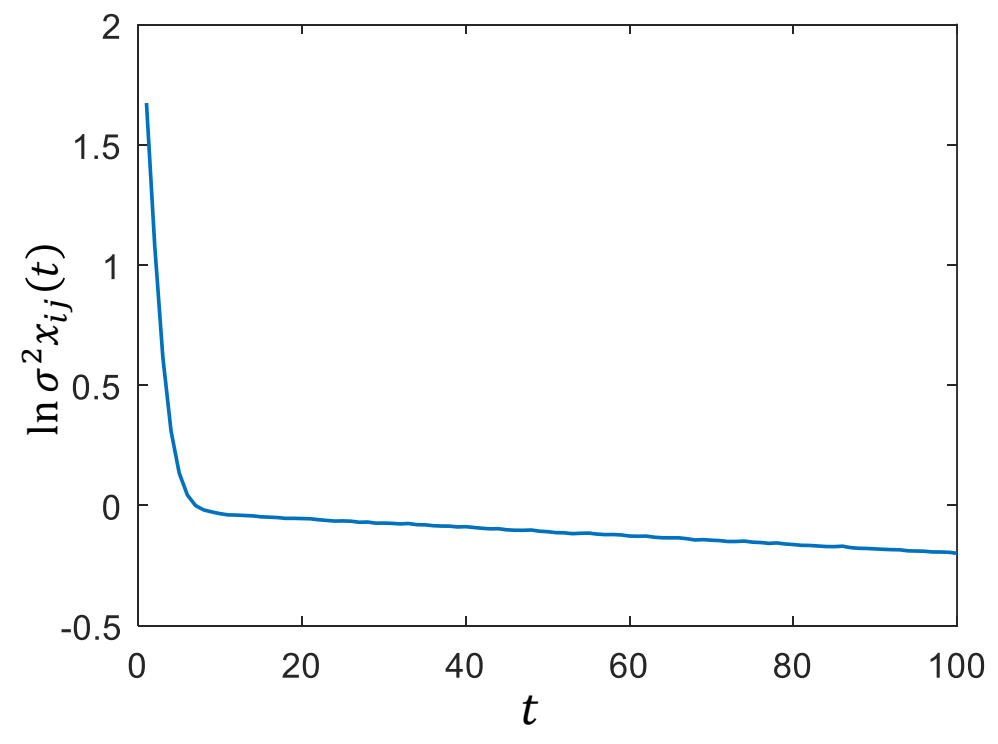

(a)

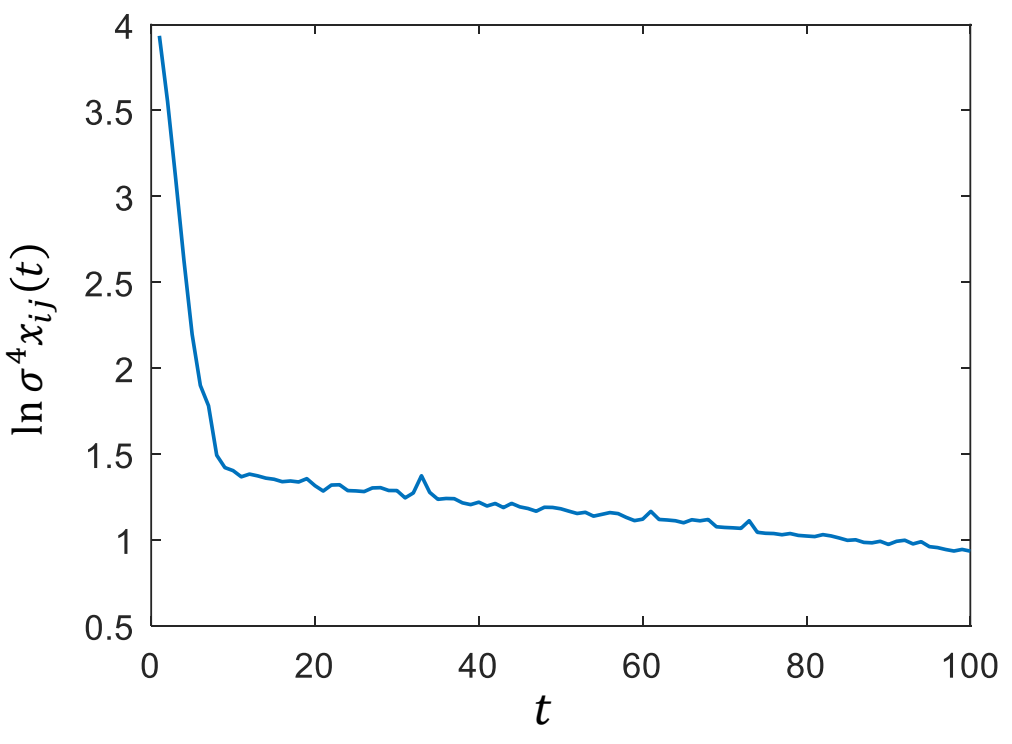

(b)

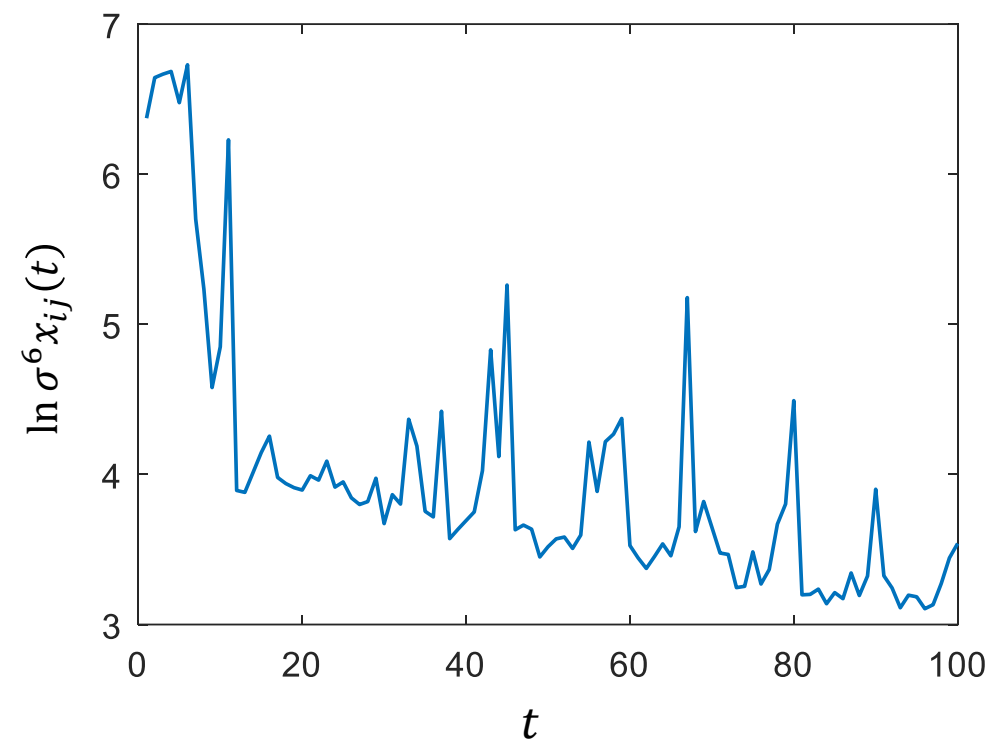

(c)

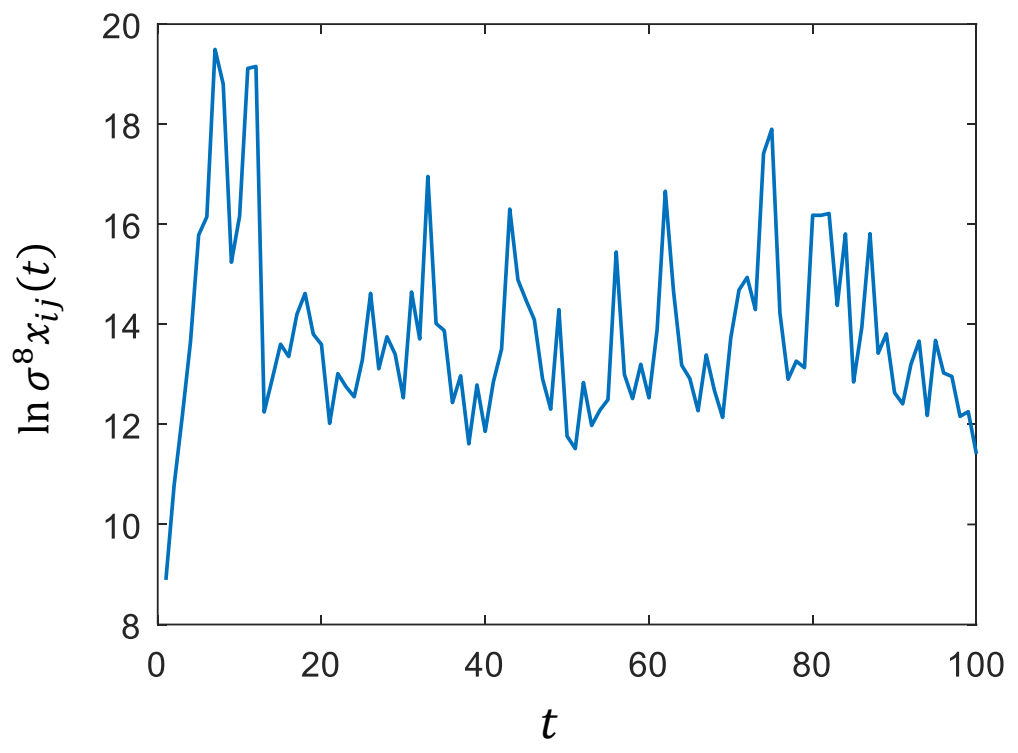

(d)

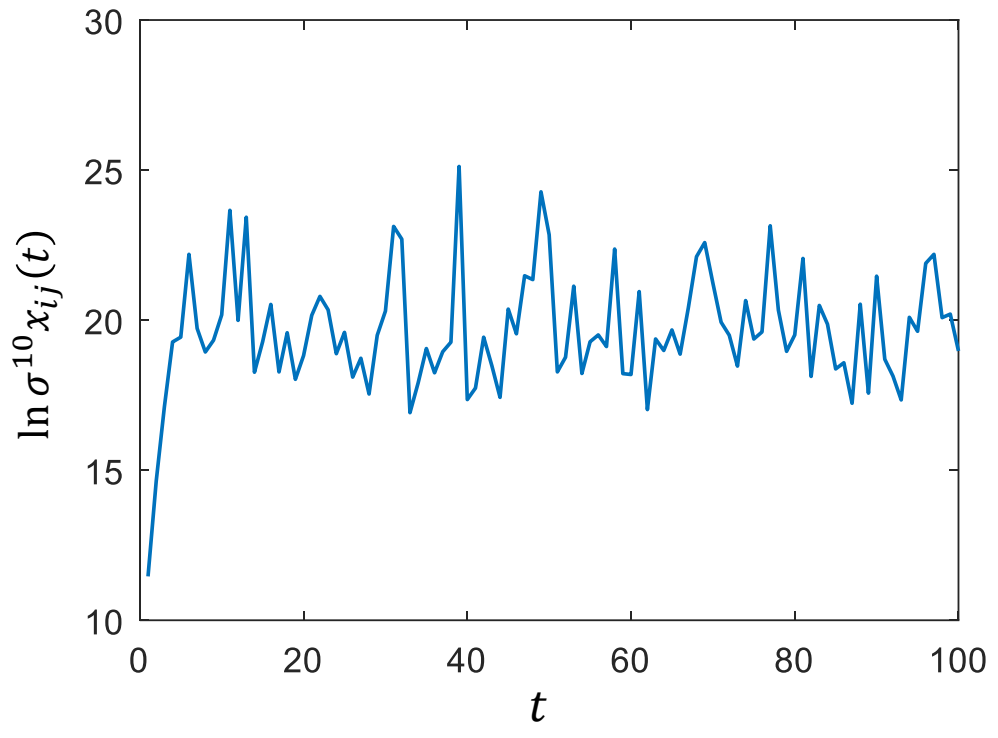

(e)

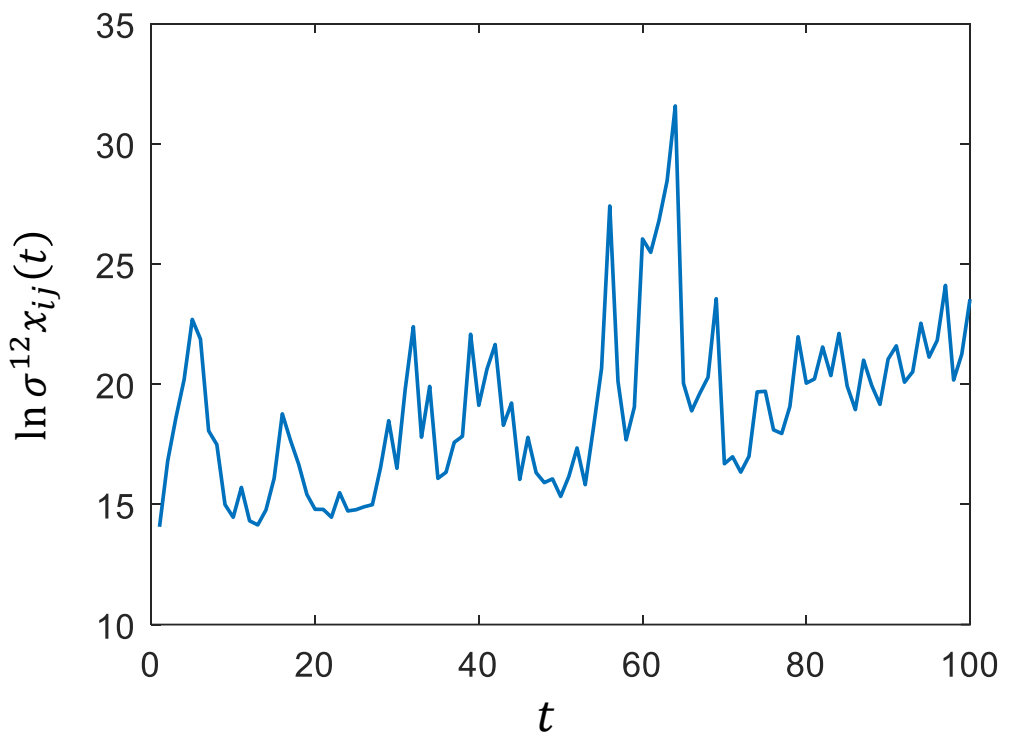

(f)

Fig. 5 $\sigma^r x_{ij}(t)$ during iterations. (a): $r = 2$, (b): $r = 4$, (c): $r = 6$, (d): $r = 8$, (e): $r = 10$, (f): $r = 12$.

It can be seen from Fig. 5 that when $r \leq 6$, the value of $\sigma^r x_{ij}(t)$ converges well, while for the cases that $r \geq 8$, the value of $\sigma^r x_{ij}(t)$ vibrates sharply and tends to increase, which are in accordance with the lemma 4. The sharp vibration of the curve of $\sigma^r x_{ij}(t)$ when $r \geq 8$ is caused by system errors and the little distance between $a$ (from 2 to 1.9), and the threshold convergence value of $\sigma^r x_{ij}(t)$ ($a_r = 1.87$ when $r = 8$, please see Fig. 3 for more details).

### 5.1.2 Verification of the lemma 6

In order to verify the lemma 6, the GWO under stagnation assumption is performed from $a = 2$ to $a = 1.9$, and the PDFs of the $x_{ij}(t)$ during iterations are simulated. Here, we set $T = 2000$, and the values of $p_{kj}$ are randomly chosen. The GWO under stagnation assumption is run for $10^6$ trials from $t = 1$ to $t = \lfloor 0.05T \rfloor = 100$ (the initial position $x_{ij}(1) \sim U[-4,4]$ in each trial), and the upper-side midpoints of all rectangles in each frequency histogram of $x_{ij}(t)(t = 2,10,20,30,50,100)$ are connected to simulate the curves of the PDFs of $h_{ij,t}(u_j)(t = 2,10,20,30,50,100)$ as shown in Fig. 6.

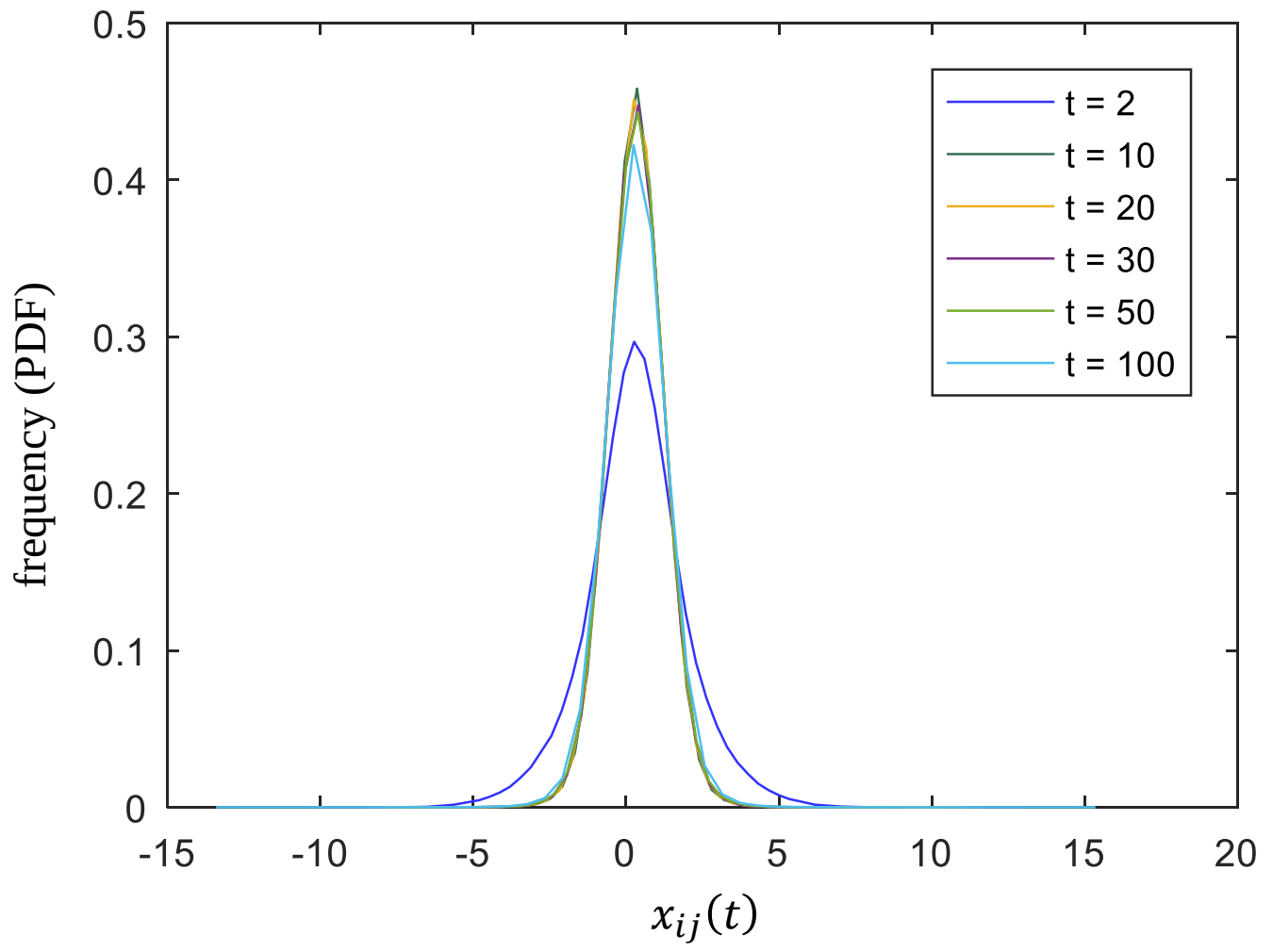

Fig. 6 PDFs of $x_{ij}(t)(t = 2,10,20,30,50,100)$

It can be seen from Fig. 6 that the shape of the curve of $h_{ij,t}(u_j)$ remains almost unchanged during iterations, which indicates that $\forall u_0$, the value of $h_{ij,\infty}(u_0)$ will not tend to 0, where $h_{ij,\infty}(u_0)$ denotes the PDF of the limit distribution of $x_{ij}(t)$ when $T \to \infty$ and $t \to \lfloor 0.05T \rfloor$. However, in Fig. 6, the domains of all $h_{ij,t}(u_j)$ are within $(-15,15)$, which is because the value of $h_{ij,t}(u_j)$ becomes too small to be simulated within $10^6$ samples when the value of $\left|u_j - \frac{1}{3}\sum p_j\right|$ is large enough (please see the Fig. 8 in the part I of the paper for further discussions of the small value of $h_{ij,t}(u_j)$ when $u_j$ is close to the upper/lower bound of its domain).

### 5.1.3 Verification of the lemma 7

In order to verify the lemma 7, the GWO under stagnation assumption is performed from $a = 2$ to $a = 1.9$, and the probabilities of finding a solution within the $j$th dimension of several given neighborhoods of the global optimum $[lb_j, ub_j]$ are simulated. Here, we set $T = 2000$, the values of $p_{kj}$ are randomly chosen such that $\frac{1}{3}\sum p_j = 0$, the interval length $|ub_j - lb_j|$ is set to be 1, and the values of $lb_j$ are set as -0.5, 0.5, 1.5, 2.5, 3.5, 4.5, 5.5, respectively. For each group of values of $ub_j$ and $lb_j$, the GWO under stagnation assumption is run for $10^6$ trials from $t = 1$ to $t = \lfloor 0.05T \rfloor = 100$, by dividing $10^6$ into the trials that $x_{ij}(t)$ falls into $[lb_j, ub_j]$ (which means the $x_{ij}(t)$ succeeds in approximately finding the global optimum, thus replaced by finding the global optimum in following discussions), the approximate value of the probabilities of finding the global optimum at iteration $t$ is obtained. Since the values of the probabilities are sometimes relatively small, the logarithm of the probability is calculated and shown in Fig. 7.

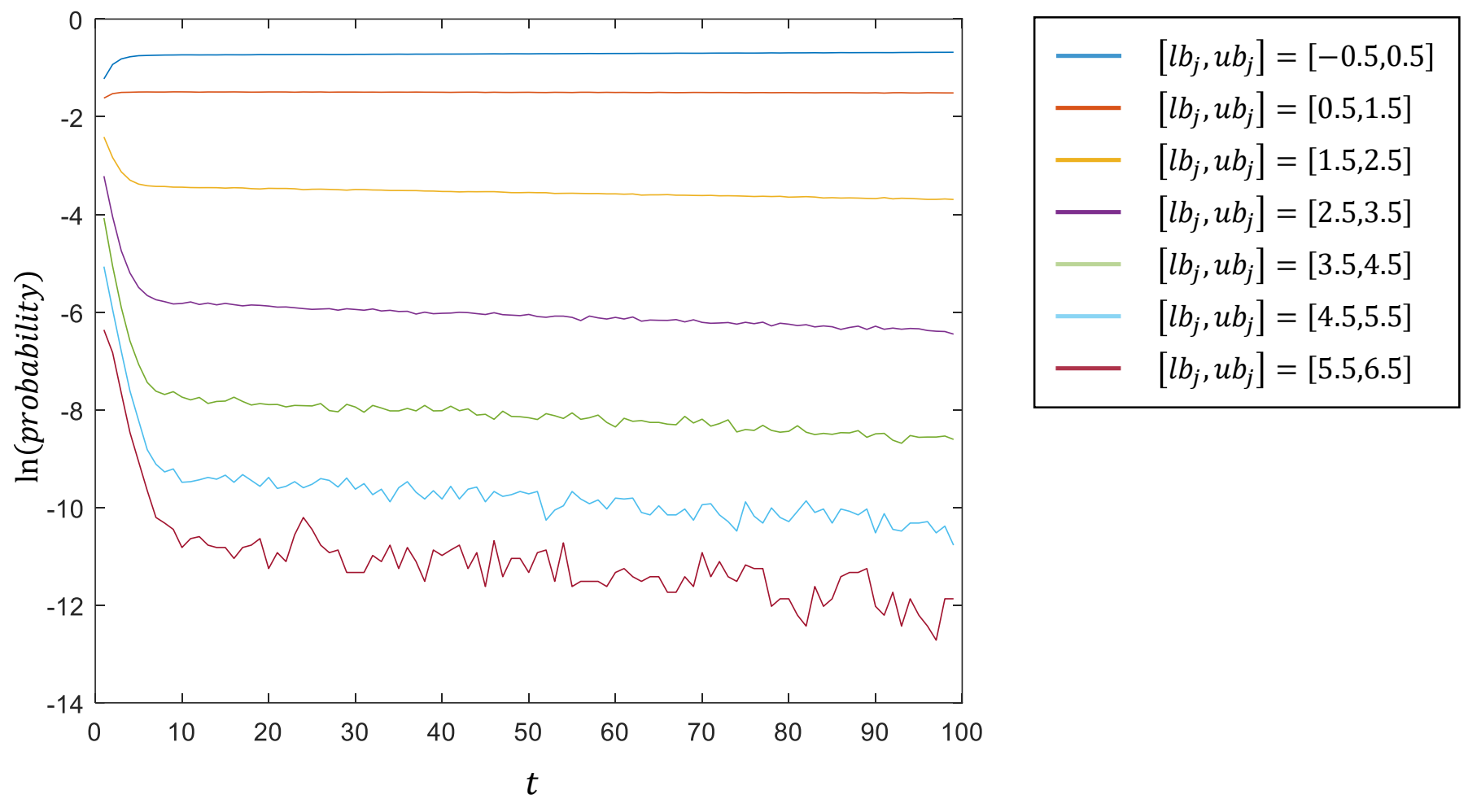

Fig. 7 Probability of finding the global optimum during iterations

Since $h_{ij,t}(u_j)$ is a single-peak function symmetrical about $\frac{1}{3}\sum p_j = 0$, obviously the probability of finding the global optimum decreases with the increase of the distance between $\frac{1}{3}\sum p_j = 0$ and the neighborhood $[lb_j, ub_j]$, which is shown in Fig. 7. In addition, it can be seen that the value of the probability of finding the global optimum slightly decreases during iterations, which is due to the slight decrease of $a$ (from 2 to 1.9). However, when $t \to \lfloor 0.05T \rfloor$, the probability tends to a finite positive number rather than 0, which is in accordance with the Lemma 8. If the value of $a$ remains unchanged during iterations, the corresponding simulation results show that the probability will converge to a positive number without the decreasing process (as shown in Fig. 7), which are shown in Fig. 8 (all simulation parameters are identical to those in obtaining Fig. 7 except for $a = 2$).

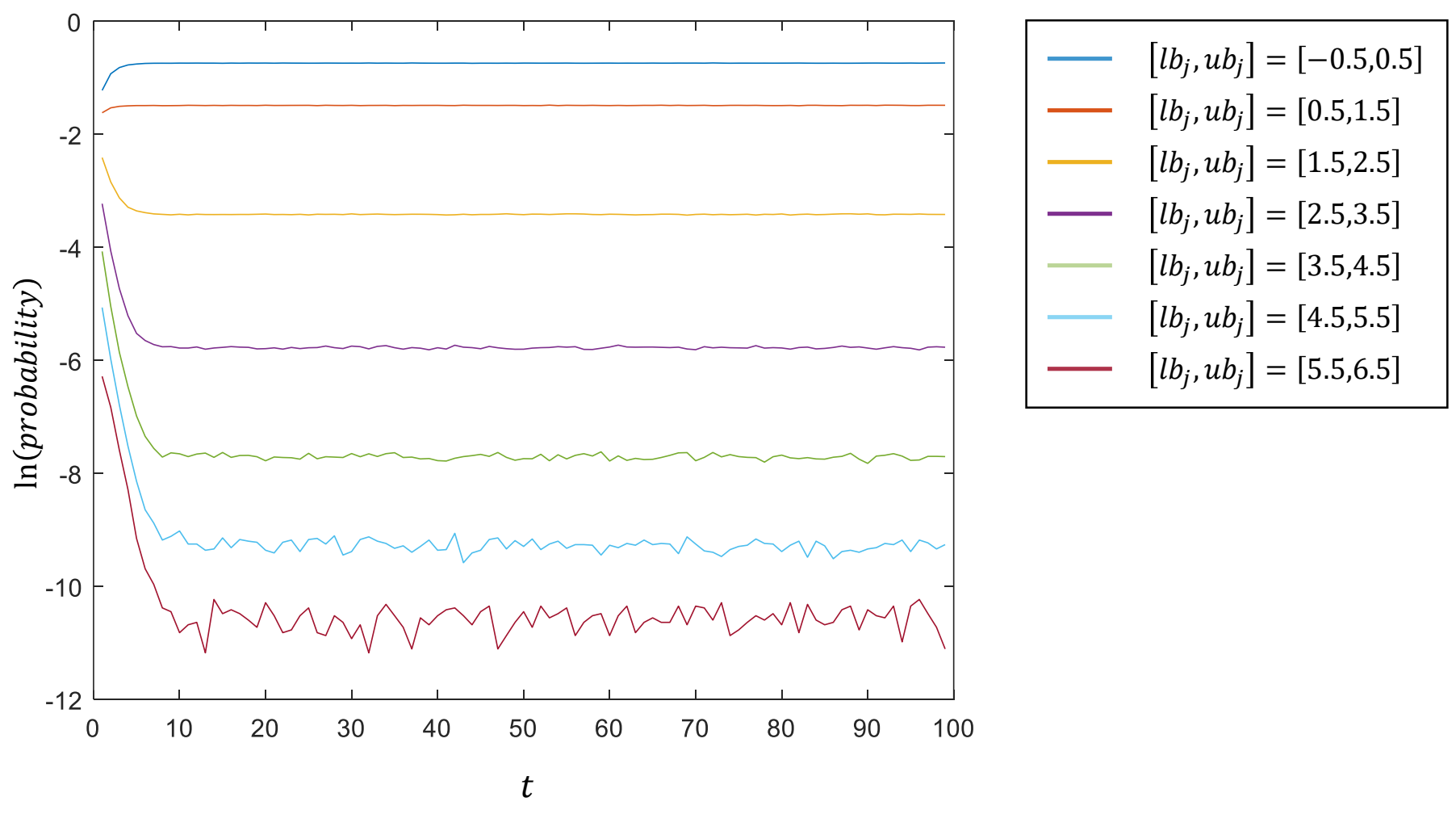

Fig. 8 Probability of finding the global optimum during iterations ($a = 2$).

It should be noticed that when the distance between the mutation center $\frac{1}{3}\sum p_j = 0$ and the neighborhood $[lb_j, ub_j]$ becomes larger than 5.5 (the largest distance as shown in Fig. 7), the probability becomes too small to be simulated within $10^6$ trials. Thus, in order to detect the probability under the cases that the mutation center is far away from the neighborhood, the GWO under stagnation assumption is run for $5 \times 10^8$ trials, the values of $lb_j$ are set as 6.5, 7.5, 8.5, 9.5, 10.5, 11.5, 12.5, 13.5, respectively, and the simulation results are shown in Fig. 9 (other simulation parameters are identical to those in obtaining Fig. 7).

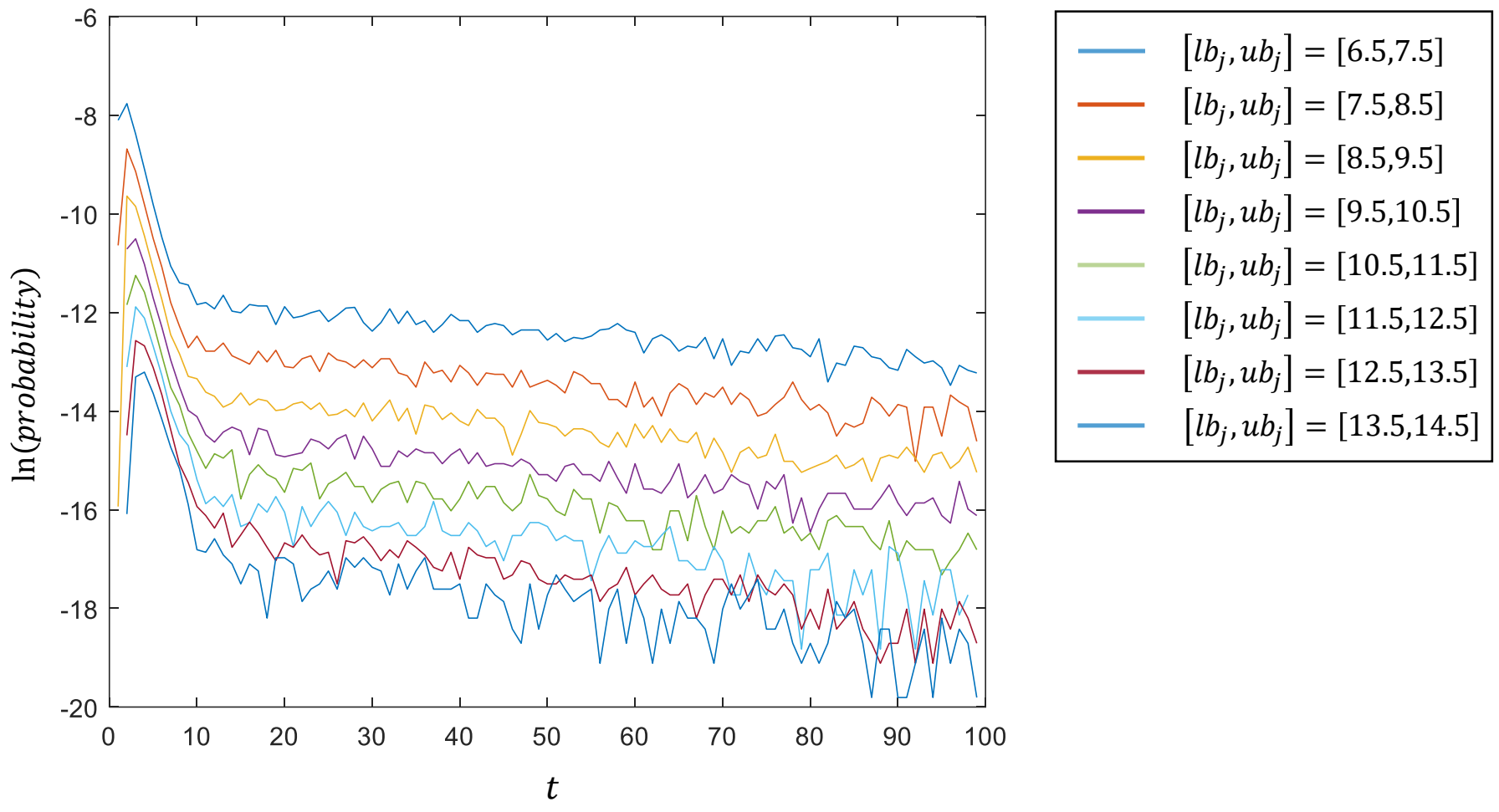

Fig. 9 Probability of finding the global optimum that is far away from the mutation center ($\frac{1}{3}\sum p_j$) during iterations.

It can be seen from Fig. 9 that although the probability decreases sharply with the increase of the distance between 0 and the neighborhood, it still converges to a positive number rather than 0. In conclusion, regardless of the distance between the stagnation point and the global optimum, the global searching ability of each search agent is always guaranteed during iterations, indicating the global convergence of the GWO under stagnation assumption.

## 5.2 Discussions

In section 4, the probability-1 global convergence of the GWO is analyzed under the stagnation assumption. Based on the stagnation assumption, the mutation center of each search agent ($\frac{1}{3}\sum p_j$) becomes constant, and the properties of the distribution of the new solution become easier to analyze. However, regarding the original GWO, the positions of the best three search agents keep changing during iterations (unless the best solution never updates, i.e. the fitness curve becomes a horizontal line, which is almost impossible in actual situations), which makes the mechanism of global searching slightly different from that of the GWO under stagnation assumption. Nevertheless, if we can clarify that the probability of finding the global optimum for the GWO under stagnation is less than that for the original GWO, it can be assumed that the probability-1 global convergence of the original GWO is guaranteed as well. It should be noted that this section is not intended to provide a rigorous theoretical proof for the probability-1 global convergence of the original GWO. More research on the theoretical analysis of the original GWO will require more in-depth exploration and exploitation in the near future.

Let's consider the positions of the best three search agents $\mathbf{p_k}(t)(k = 1,2,3)$. In the original GWO, the new solutions are generated around $\mathbf{p_k}(t)$. If a better solution is found, then the positions of $\mathbf{p_k}(t)$ will be updated, and the agents will move around the new best search agents $\mathbf{p_k}(t)$. Obviously, the updating mechanism of $\mathbf{p_k}(t)$ guarantees that the whole swarm is moving towards the global optimum $\mathbf{x_g}$ under the guidance of $\mathbf{p_k}(t)$, which moves towards the $\mathbf{x_g}$ itself. Such mechanism substantially accelerates the global optimum searching process. While in the GWO under stagnation assumption, the positions of $\mathbf{p_k}(t)$ remain unchanged (if the global optimum is not found), then the process that the swarm moves towards $\mathbf{x_g}$ under the guidance of $\mathbf{p_k}(t)$ no longer exists, which means the corresponding global optimum searching ability of the population is substantially reduced.

To verify the aforementioned conclusions, the corresponding numerical simulations are conducted in the last part of this section. During simulations, ten classical test problems that are widely utilized by many researchers [38-39] with box feasible regions are optimized by both the original GWO and the GWO under stagnation assumption respectively, containing both single-peak functions and multi-peak functions with different characteristics. The information of these test problems are given in Table III. Here, the number of search agents $N = 10$, the total number of iterations $T = 50000$, and the initial positions of search agents $\mathbf{x_i}(1)$ are randomly chosen within the feasible region. The neighborhood of the $\mathbf{x_g}$ is defined as $S_\varepsilon = \left\{\mathbf{x_0} \middle| \frac{f(\mathbf{x_0})-f(\mathbf{x_g})}{\min_i f(\mathbf{x_i}(1))-f(\mathbf{x_g})} \le \varepsilon, 0 < \varepsilon < 1\right\}$, and the optimization processes of the original GWO and GWO under stagnation assumption are terminated if a feasible solution $\mathbf{x_0} \in S_\varepsilon$ is found, in which case we assume that the global optimum is found. For each case of the value of $\varepsilon$ (to be set as 0.1 or 0.01 respectively), each test problem is solved independently 50 trials by the original GWO and the GWO under stagnation assumption, respectively. The number of times that the original GWO or the GWO under stagnation assumption succeeds in finding the global optimum of each test problem and the corresponding average running time are provided in Table IV, respectively. In addition, in order to better illustrate the population dynamics in the original GWO, both the fitness value of the best search agent and the distance between the best search agent and the $\mathbf{x_g}$ during iterations in one trial of the original GWO before finding a solution within $S_{\varepsilon=0.001}$ or during 50000 iterations are shown in Fig. 10 (a)-(j), respectively. The dynamics of the position of the best search agent during iterations of the GWO under stagnation assumption is not provided since it remains unchanged during iterations before finding the global optimum.

Table III Benchmark functions

| Test problem No. | Objective function | $D$ | Box constraints | $\mathbf{x_g}$ | $f(\mathbf{x_g})$ | Problem type* |
|---|---|---|---|---|---|---|
| 1 | $f_1(\mathbf{x}) = \sum_{i=1}^{D} x_i^2$ | 3 | $[-3,3]^D$ | $[0,\ldots,0]$ | 0 | S |
| 2 | $f_2(\mathbf{x}) = \sqrt{\sum_{i=1}^{D} \lvert x_i \rvert^{2+4\frac{i-1}{D-1}}}$ | 3 | $[-3,3]^D$ | $[0,\ldots,0]$ | 0 | S |
| 3 | $f_3(\mathbf{x}) = \sum_{i=1}^{D} \lvert x_i \rvert + \prod_{i=1}^{D} \lvert x_i \rvert$ | 3 | $[-3,3]^D$ | $[0,\ldots,0]$ | 0 | S |
| 4 | $f_4(\mathbf{x}) = \sum_{i=1}^{D-1} [100(x_i^2 - x_{i+1})^2 + (1-x_i)^2]$ | 3 | $[-3,3]^D$ | $[1,\ldots,1]$ | 0 | S |
| 5 | $f_5(\mathbf{x}) = \frac{1}{4000}\sum_{i=1}^{D} x_i^2 - \prod_{i=1}^{D} \cos\frac{x_i}{\sqrt{i}} + 1$ | 3 | $[-10,10]^D$ | $[0,\ldots,0]$ | 0 | M |
| 6 | $f_6(\mathbf{x}) = \sum_{i=1}^{D} [x_i^2 - 10\cos 2\pi x_i + 10]$ | 3 | $[-10,10]^D$ | $[0,\ldots,0]$ | 0 | M |
| 7 | $f_7(\mathbf{x}) = -\sum_{i=1}^{D} x_i \sin\sqrt{\lvert x_i \rvert}$ | 3 | $[-500,500]^D$ | $420.97[1,\ldots,1]$ | $-418.98D$ | M |
| 8 | $f_8(\mathbf{x}) = \left\lvert \sum_{i=1}^{D} x_i^2 - D \right\rvert^{\frac{1}{4}} + \frac{1}{D}\left(0.5\sum_{i=1}^{D} x_i^2 + \sum_{i=1}^{D} x_i\right) + 0.5$ | 3 | $[-3,3]^D$ | $-[1,\ldots,1]$ | 0 | M |
| 9 | $f_9(\mathbf{x}) = \left(x_2 - \frac{5.1}{4\pi^2}x_1^2 + \frac{5}{\pi}x_1 - 6\right)^2 + 10\left(1 - \frac{1}{8\pi}\right)\cos x_1 + 10$ | 2 | $[-20,20]^2$ | $x_1 = 2k\pi + \pi, x_2 = \frac{5.1}{4\pi^2}x_1^2 - \frac{5}{\pi}x_1 + 6$ | 0.3979 | M |
| 10 | $f_{10}(\mathbf{x}) = \left(4 - 2.1x_1^2 + \frac{1}{3}x_1^4\right)x_1^2 + x_1x_2 + (-4 + 4x_2^2)x_2^2$ | 2 | $[-3,3]^2$ | $[\pm 0.09, \mp 0.71]$ | -1.0316 | M |

* S denotes single-peak characteristic, M denotes multi-peak characteristic.

Table IV Simulation results

| Test problem No. | The number of times succeeding in finding the global optimum | | | | Average number of iterations* | | | |
|---|---|---|---|---|---|---|---|---|
| | GWO | | GWO under stagnation assumption | | GWO | | GWO under stagnation assumption | |
| | $\sigma = 0.1$ | $\sigma = 0.01$ | $\sigma = 0.1$ | $\sigma = 0.01$ | $\sigma = 0.1$ | $\sigma = 0.01$ | $\sigma = 0.1$ | $\sigma = 0.01$** |
| 1 | 50 | 50 | 50 | 48 | 4 | 7 | 44 | 998 |
| 2 | 50 | 50 | 49 | 46 | 5 | 9 | 232 | 2243 |
| 3 | 50 | 50 | 49 | 10 | 7 | 13 | 491 | – |
| 4 | 50 | 50 | 49 | 44 | 10 | 67 | 240 | 1738 |
| 5 | 50 | 50 | 49 | 20 | 180 | 2240 | 2058 | 33250 |
| 6 | 50 | 50 | 50 | 9 | 18 | 55 | 970 | – |
| 7 | 12 | 4 | 0 | 0 | 38316 | – | – | – |
| 8 | 45 | 3 | 0 | 0 | 7967 | – | – | – |
| 9 | 50 | 50 | 50 | 50 | 13 | 83 | 17 | 109 |
| 10 | 50 | 50 | 50 | 50 | 4 | 12 | 30 | 148 |

* If the solution that satisfies Eq.(3.2) is not found after 50000 iterations, the trial will be eliminated

** If the number of times that the GWO under stagnation assumption succeeds in finding the global optimum is less than 10 times in a total of 50 trials, it can be deduced that the GWO under stagnation assumption is unlikely to find the global optimum within 50000 iterations

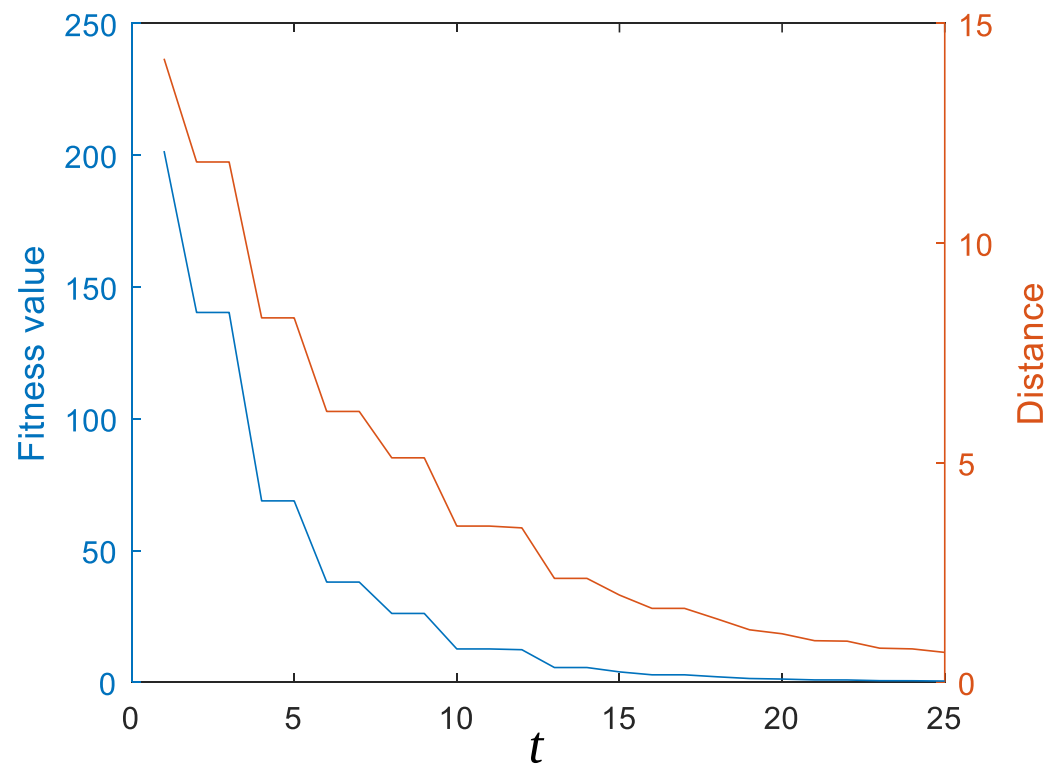

(a)

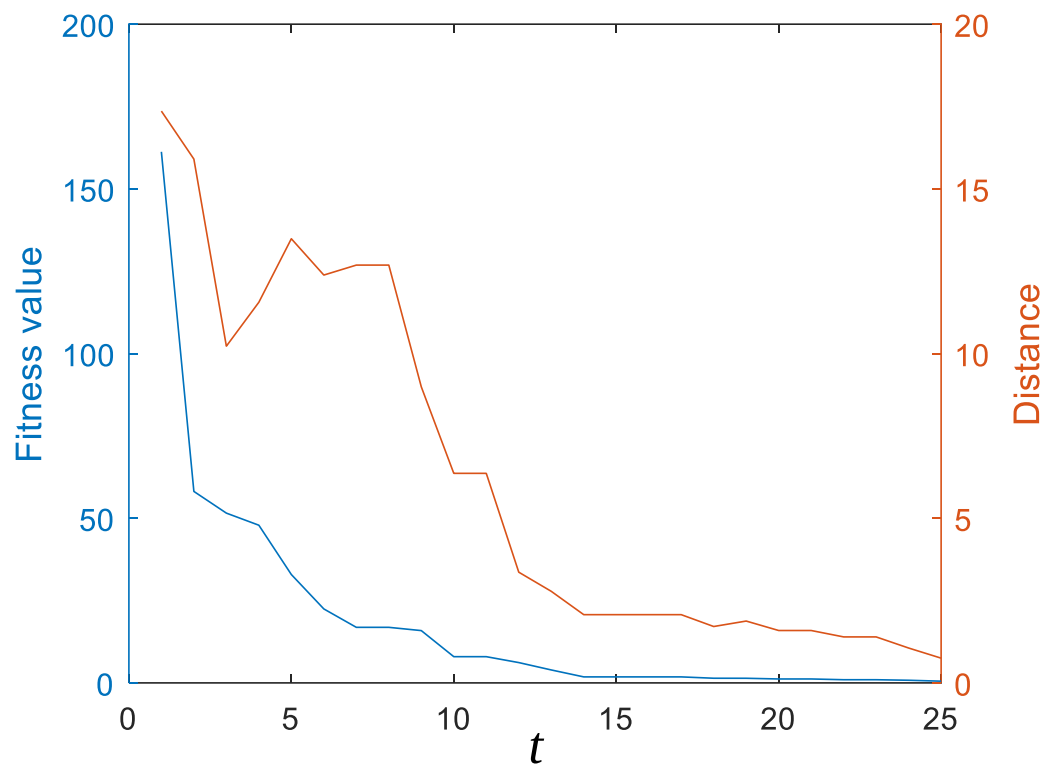

(b)

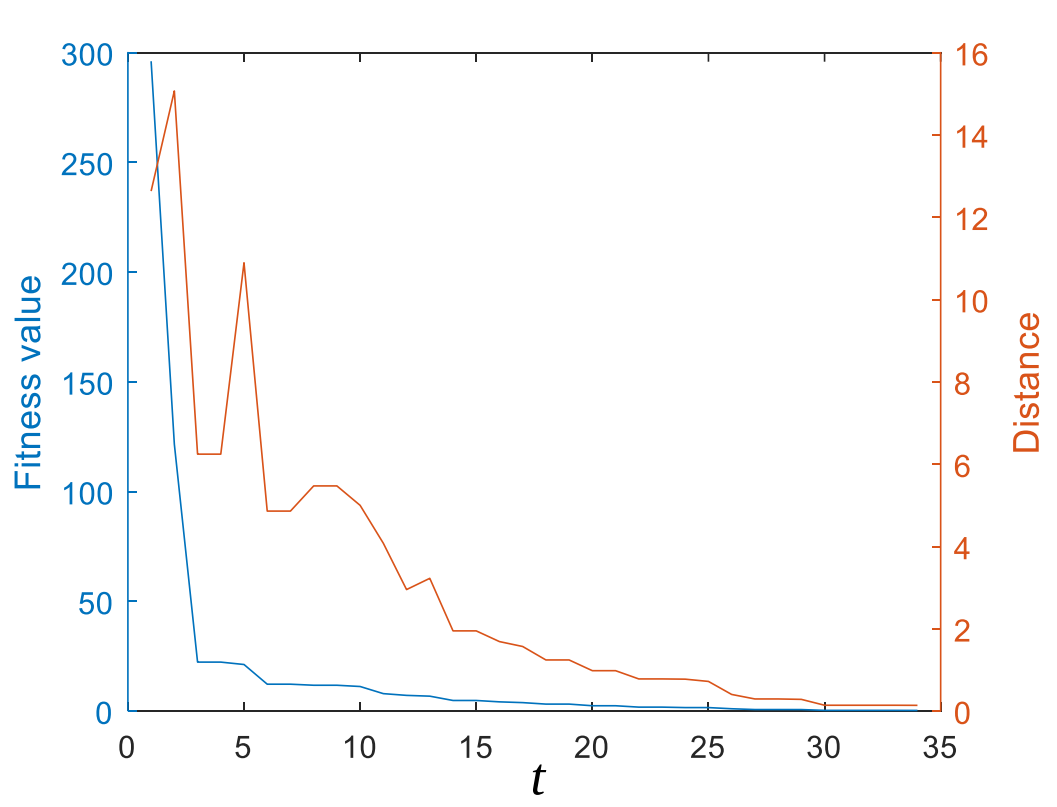

(c)

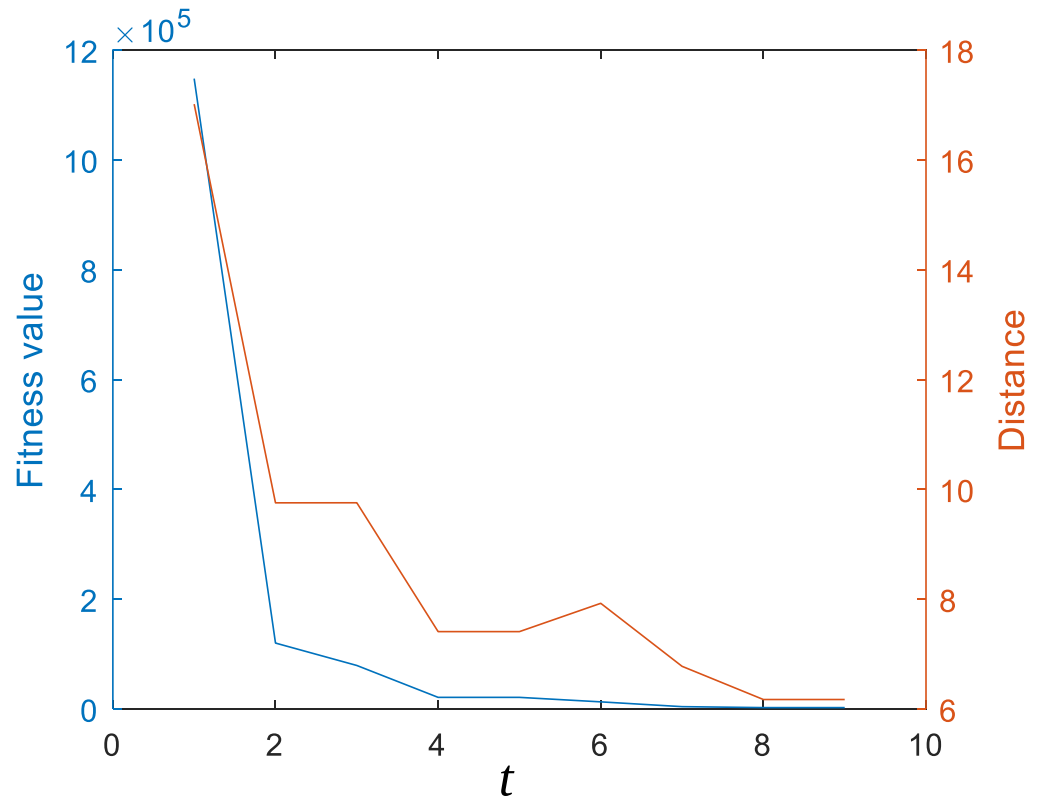

(d)

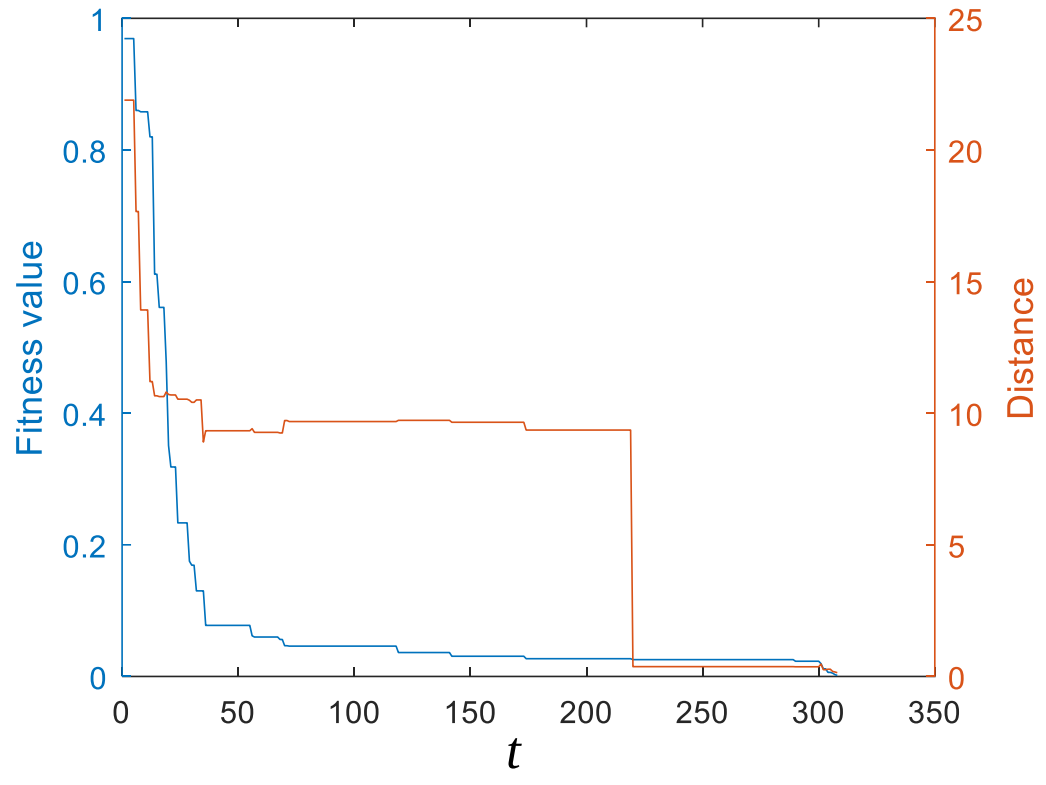

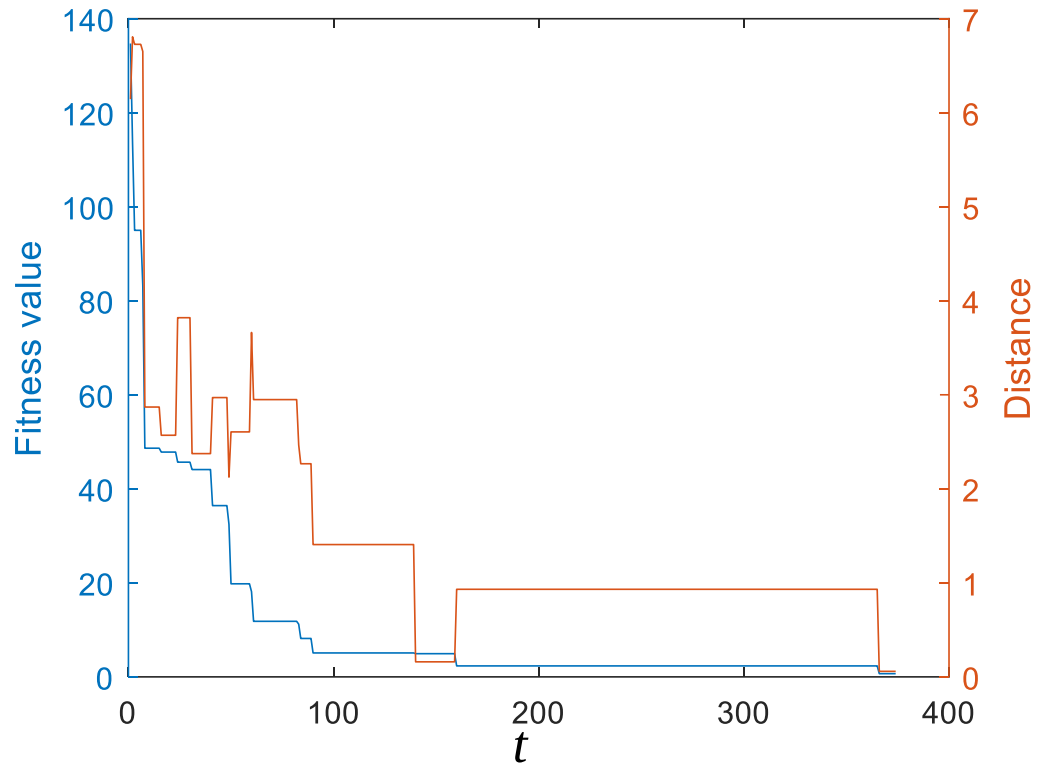

(e) (f)

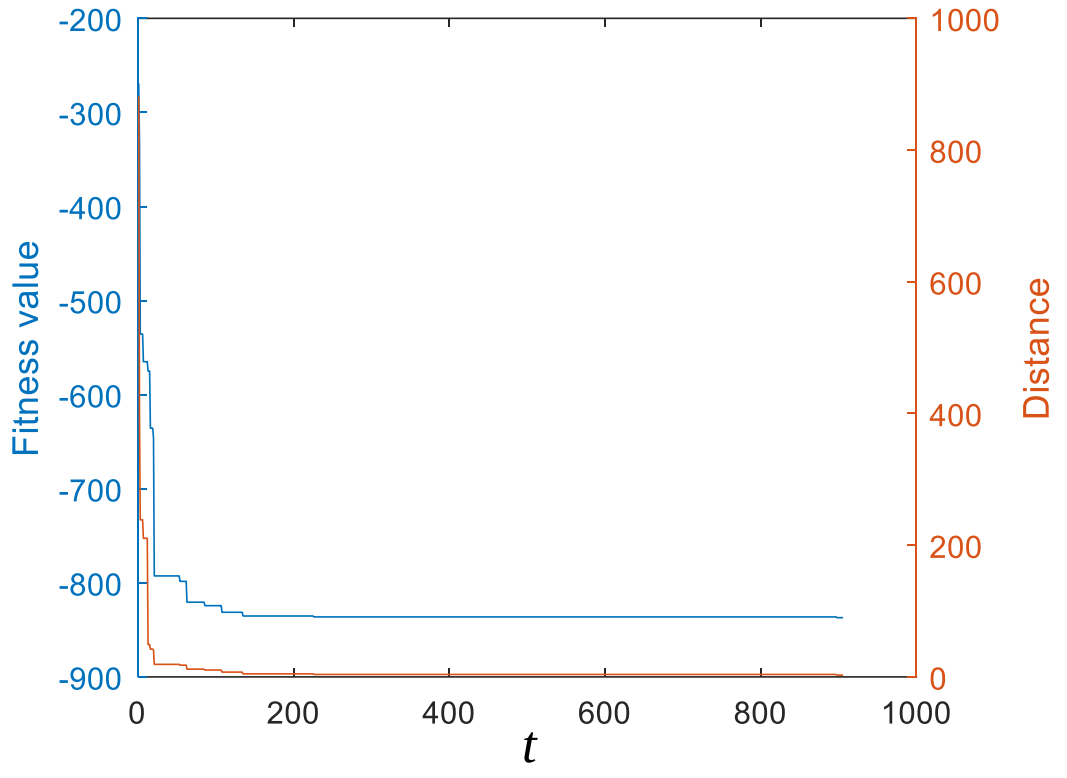

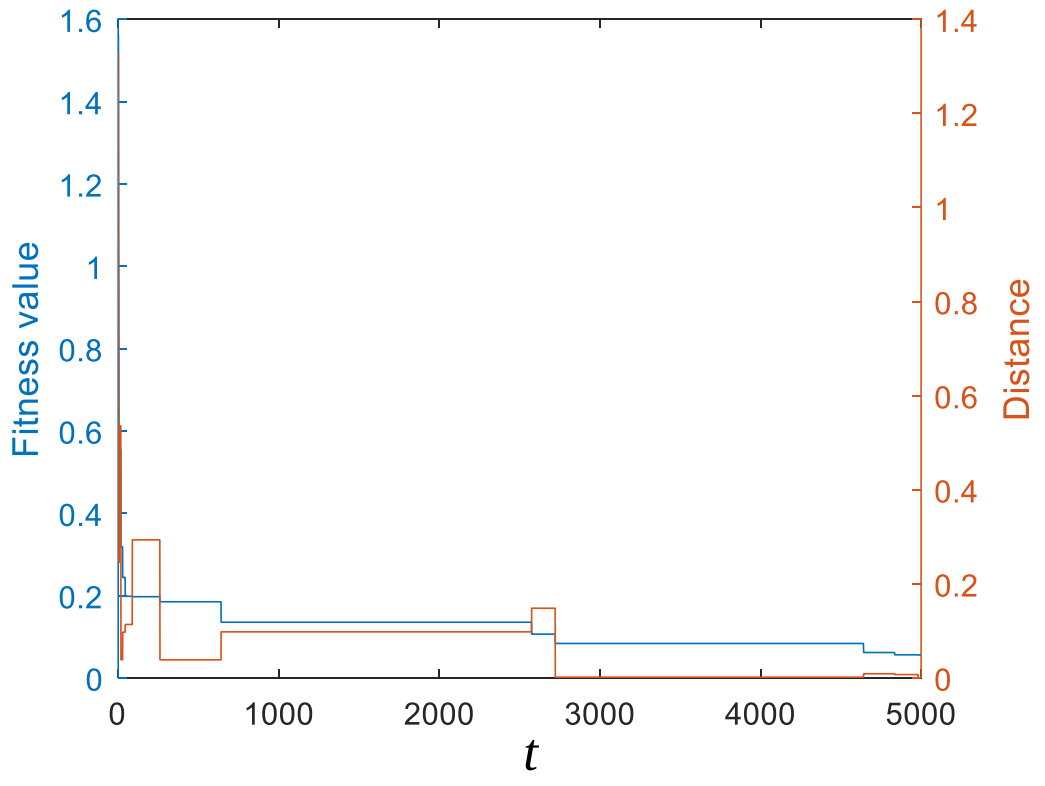

(g) (h)

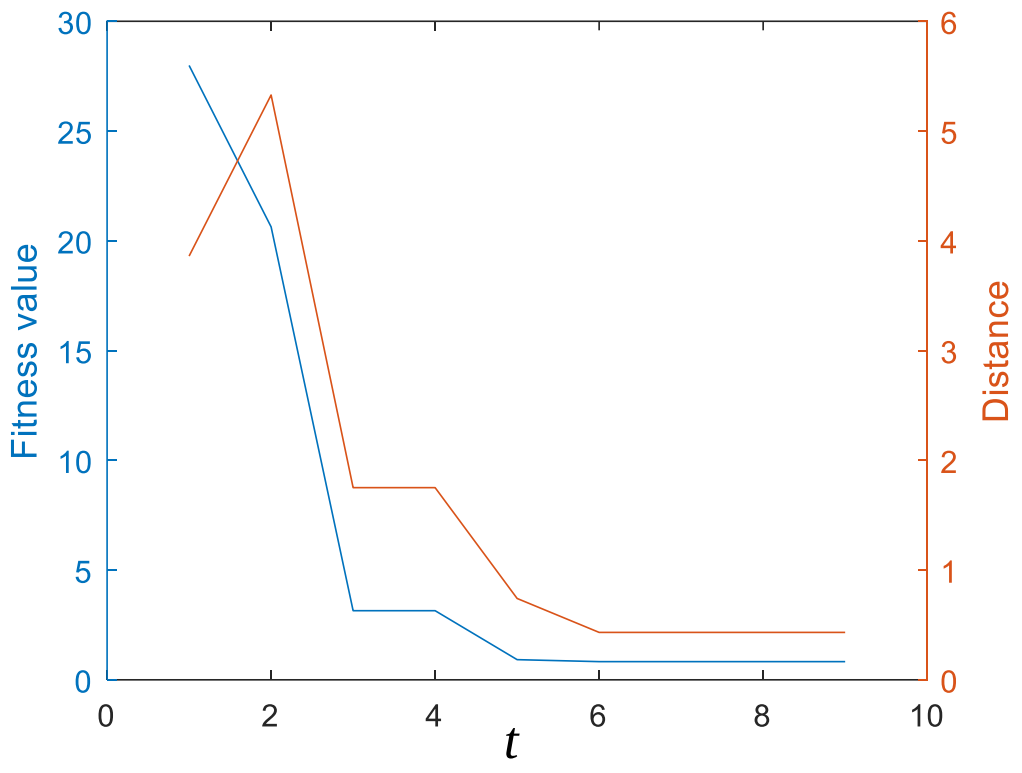

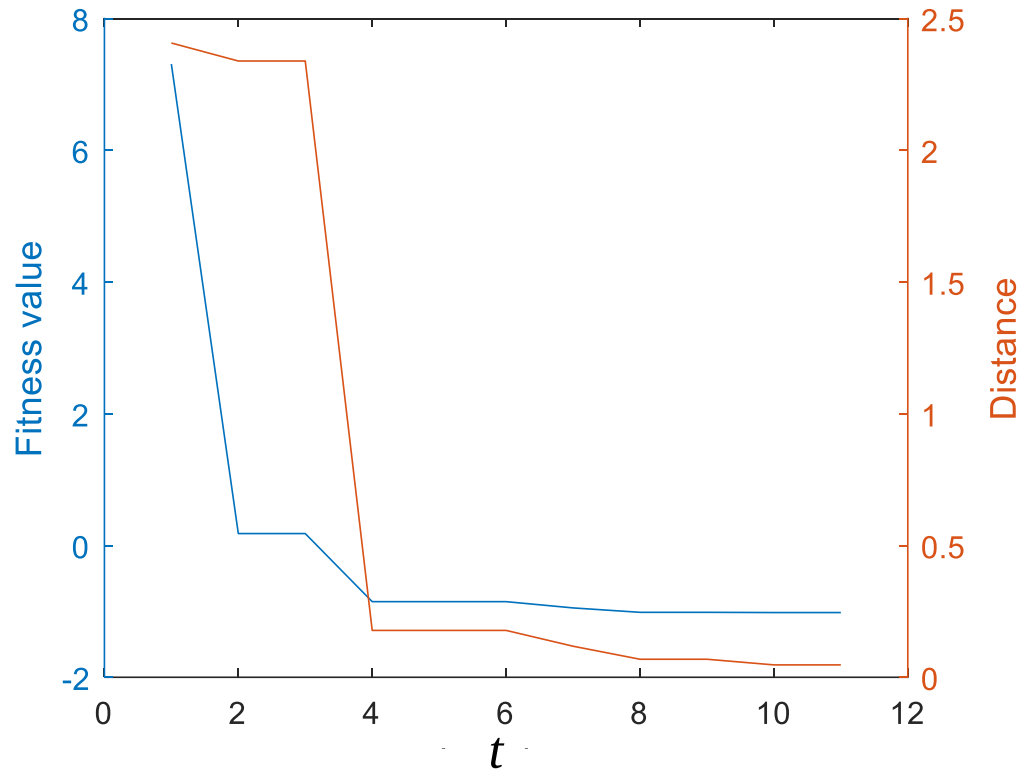

(i) (j)

Fig. 10 Fitness value of the best search agent and the distance between the best search agent and $\mathbf{x_g}$ in the GWO, (a)-(j) correspond to test problems No. 1-No. 10.

It can be seen from Table IV that for fixed $\varepsilon$ and problems that are easy to solve (such as problems No. 1 and No. 2), the success rate of approximately finding the global optimum before 50000 iterations is almost 100% for the original GWO, while there exist cases that the GWO under stagnation assumption fails to find the global optimum within 50000 iterations. As for the problems that are difficult to solve (such as problems No. 7 - No. 8), the success rate of approximately finding the global optimum before 50000 iterations is almost 0% for the GWO under stagnation assumption. In addition, the average number of iterations before finding the global optimum for the original GWO is much less than that of the GWO under stagnation assumption. In conclusion, regardless of the value of $\varepsilon$ and the number of test problems, the original GWO is superior to the GWO under stagnation assumption. In addition, it can be seen from Fig. 10 that as the number of iterations increases and the fitness value of the best search agent decreases, the distance between the best search agent and $\mathbf{x_g}$ tends to decrease (the reason that the fitness value is not a strictly decreasing function of $t$ is that for multi-modal problems, the best solution is likely to temporarily fall into the local optimum) corresponding to the conclusion mentioned above. All these simulation results show that the global searching ability of the original GWO is much better than that of the GWO under stagnation assumption, indicating that the global convergence of the original GWO is guaranteed as well.

# 6. Conclusion

In the part II of the paper, based on the theoretical findings discussed in the part I such as the sampling distribution and probabilistic stability of the GWO, we analyze the global convergence properties of the GWO. The stagnation assumption that the positions of the best three search agents $\mathbf{p_k}(t)(k = 1,2,3)$ remain unchanged during iterations runs throughout the entire proof process.

In this study, the convergence properties of the GWO under stagnation assumption are abstracted and modelled into two propositions, corresponding to *global searching ability analysis* and *probability-1 global convergence analysis*, respectively. The global searching ability of the GWO under stagnation assumption is firstly proved and carefully studied.

In order to carry out the probability-1 global convergence analysis of the GWO under stagnation assumption, two concepts called algorithm condition and convergence condition are introduced and studied. The characteristics of the distribution of new solution in the GWO under stagnation assumption is further discussed, and the conclusion that the $r$-order central moment of $x_{ij}(t)$ tends to infinity when $t$ tends to infinity for any positive even number $r \geq 8$ is proved and serves as a crucial lemma. Based on this lemma, the GWO under stagnation assumption is proved to satisfy both algorithm condition and convergence condition.

All the results proved above are verified by numerical simulations. In addition, although the convergence analysis is carried out under the stagnation assumption, it can be seen from qualitative analysis and numerical simulations that the global convergence property of the original GWO is superior to that of the GWO under stagnation, thus the global convergence of the original GWO can be asserted as well, which will be proved strictly in our future work.

# Appendix

## Proof of proposition 1

Mark the widest searching scale of $x_{ij}(t)$ as $d_j(t) = \max_{x_{ij}(t)} \left| x_{ij}(t) - \frac{1}{3}\sum p_j \right|$ (i.e. the domain

of $h_{ij,t}(u_j)$ becomes $\left(\frac{1}{3}\sum p_j - d_j(t), \frac{1}{3}\sum p_j + d_j(t)\right)$), and it has been proved in the part I (please see the Proposition 3.1 in the part I of the paper) that when $t \geq 2$, $d_j(t) \geq d_j(1)\prod_{\tau=1}^{t-1} 2\left(1-\frac{\tau}{T}\right)$ holds. Mark the $j$th component of the maximum distance between $\frac{1}{3}\sum_{k=1}^{3}\mathbf{p_k}$ and $S_0$ as $d_{0j} = \max\left\{\left|lb_j - \frac{1}{3}\sum p_j\right|, \left|ub_j - \frac{1}{3}\sum p_j\right|\right\}$, suppose that $T > 4t$ (i.e. $a > 1.5$), then $d_j(t) \geq 1.5^{t-1}d_j(1)$. Set $t_0 = \left\lceil \log_{1.5}\frac{d_{0j}}{d_j(1)} + 1\right\rceil, T > 4t_0$, then it can be deduced that $d_j(t_0) \geq 1.5^{t_0-1}d_j(1) \geq 1.5^{\log_{1.5}\frac{d_{0j}}{d_j(1)}}d_j(1) = d_{0j}$, accordingly $\forall t = t_0, t_0+1, \dots, \left\lfloor\frac{T}{2}\right\rfloor$, the domain of $x_{ij}(t)$ contains $[lb_j, ub_j]$, in which case $\mathbb{P}\{x_{ij}(t) \in [lb_j, ub_j]\} = \int_{lb_j}^{ub_j} h_{ij,t}(u_j)\mathrm{d}u_j \geq (ub_j - lb_j)\min\{h_{ij,t}(lb_j), h_{ij,t}(ub_j)\} > 0$.

In conclusion, set $t_0 = \max_{j=1,\dots,d}\left\lceil \log_{1.5}\frac{d_{0j}}{d_j(1)} + 1\right\rceil, T = 4t_0 + 1$, then $\forall t = t_0, t_0 + 1, \dots, \left\lfloor\frac{T}{2}\right\rfloor, \mathbb{P}\{\mathbf{x_i}(t) \in S_0\} > 0$ holds. □

## Proof of lemma 3

(1) $f_r(t) = \frac{a^r}{3^r}\sum_{p=0}^{\frac{r}{2}}\frac{\binom{r}{2p}}{(2p+1)}\sum_{q=0}^{\frac{r}{2}-p}\frac{\binom{2\left(\frac{r}{2}-p\right)}{2q}}{(2q+1)\left[2\left(\frac{r}{2}-p\right)-2q+1\right]}$

$$= \frac{a^r}{3^r}\sum_{p=0}^{\frac{r}{2}}\frac{\binom{r}{2p}}{(2p+1)\left[2\left(\frac{r}{2}-p\right)+1\right]\left[2\left(\frac{r}{2}-p\right)+2\right]}\sum_{q=0}^{\frac{r}{2}-p}\binom{2\left(\frac{r}{2}-p\right)+2}{2q+1}$$

$$= \frac{a^r}{3^r}\sum_{p=0}^{\frac{r}{2}}\frac{2^{2\left(\frac{r}{2}-p\right)+1}\binom{r}{2p}}{(2p+1)\left[2\left(\frac{r}{2}-p\right)+1\right]\left[2\left(\frac{r}{2}-p\right)+2\right]}$$

$$= \frac{a^r}{3^r}\frac{1}{(r+1)(r+2)(r+3)}\sum_{p=0}^{\frac{r}{2}}2^{2p+1}\binom{r+3}{2p+2}$$

$$= \frac{a^r}{3^r}\frac{1}{(r+1)(r+2)(r+3)}\left(\frac{1}{2}\sum_{p=0}^{\frac{r}{2}+1}2^{2p}\binom{r+3}{2p} - 1\right)$$

$$= \frac{a^r}{3^r}\frac{3^{r+3}-3}{4(r+1)(r+2)(r+3)}$$

(2) Obviously $a_r = 3\left(\frac{4(r+1)(r+2)(r+3)}{3^{r+3}-3}\right)^{\frac{1}{r}}$. Since $\ln a_r = \frac{1}{r}\ln\frac{4(r+1)(r+2)(r+3)}{3^{r+3}-3} + \ln 3$, set $u = r+2, A = 3(3^u - 1), B = 4u(u-1)(u+1)$, we only need to prove that $g(u) = \frac{1}{u-2}\ln\frac{A}{B}$ is an increasing function of $u \geq 4$ (please notice that $r$ is an even integer). It can be calculated that $g'(u) = \frac{\frac{A'B-AB'}{AB}(u-2)-\ln\frac{A}{B}}{(u-2)^2}$. Let $h(u) = \frac{A'B-AB'}{AB}(u-2) - \ln\frac{A}{B}$, then

$$h'(u) = \frac{1}{A^2B^2}[AB(A''B(u-2) + A'B - AB''(u-2) - AB') - (A'B - AB')(u-2)(A'B + AB')] - \frac{A'B-AB'}{AB}$$

$$= \frac{u-2}{A^2B^2}\left(AA''B^2 - BB''A^2 + A^2B'^2 - A'^2B^2\right)$$

$$= \frac{144(u-2)}{A^2B^2}[(3^u-1)^2(3u^4+1) - 3^u(\ln 3)^2u^2(u^2-1)^2]$$

$$> \frac{144(u-2)}{A^2B^2}[3u^4(3^u-1)^2 - 3^u \cdot 3u^6]$$

$$= \frac{432u^4(u-2)}{A^2B^2}[(3^u-1)^2 - u^2 3^u]$$

$$> \frac{432u^4 3^u(u-2)(3^u-u^2-2)}{A^2B^2} > 0$$

Thus $h(u) > h(4) = \frac{243\ln 3-188}{120} > 0$, indicating that $g(u)$ is an increasing function. Hence, the value of $a_r$ decreases with the increase of positive even number *r*.

(3) According to (2), when $r \leq 6, a_r \geq a_6 \approx 2.05 > 2$, since $f_r(t)$ is an increasing function of $a$, $f_r(t) < 1$ holds when $a \in [1.9,2]$. When $r \geq 8, a_r \leq a_8 \approx 2.05 > 2$, thus $f_r(t) > 1$ holds when $a \in [1.9,2]$. □

## Proof of lemma 4

According to lemma 3, when $a \in [1.9,2], f_2(t), f_4(t), f_6(t) < 1, f_8(t), f_{10}(t), \ldots > 1$.

Since $\sigma^0 x_{ij}(t) = 1$, and for any odd number $r$, $\sigma^r x_{ij}(t) = 0$ holds, when $T \to \infty$ and $t \leq \lfloor 0.05T \rfloor, a \in [1.9,2]$ holds, it can be obtained that $\lim_{T\to\infty}\lim_{t\to\lfloor 0.05T\rfloor} g_2(t) < \infty$ exists, accordingly $\lim_{T\to\infty}\lim_{t\to\lfloor 0.05T\rfloor} \sigma^2 x_{ij}(t) < \infty$ exists, similarly for the values of $\sigma^4 x_{ij}(t)$ and $\sigma^6 x_{ij}(t)$. However, $f_8(t) > 1$, indicating that $\lim_{T\to\infty}\lim_{t\to\lfloor 0.05T\rfloor} \sigma^8 x_{ij}(t) = +\infty$, similarly for $\sigma^{10}x_{ij}(t), \sigma^{12}x_{ij}(t), \ldots$. □

## Proof of lemma 5

According to lemma 7, there exist a positive number $\varepsilon_{ij}$ and an infinite subsequence $\{t_l\}, l = 1,2,\ldots$, such that $\forall l, \mathbb{P}\{x_{ij}(t_l+1) \in [lb_j, ub_j]\} \geq \varepsilon_{ij}$, thus

$$\prod_{t=1}^{\infty}\prod_{i=1}^{N}\left(1 - \prod_{j=1}^{D}\mathbb{P}\{x_{ij}(t+1) \in [lb_j, ub_j]\}\right)$$
$$\leq \prod_{l=1}^{\infty}\prod_{i=1}^{N}\left(1 - \prod_{j=1}^{D}\mathbb{P}\{x_{ij}(t_l+1) \in [lb_j, ub_j]\}\right)$$
$$= \exp\sum_{l=1}^{\infty}\sum_{i=1}^{N}\ln\left(1 - \prod_{j=1}^{D}\mathbb{P}\{x_{ij}(t_l+1) \in [lb_j, ub_j]\}\right)$$
$$\leq \exp\left(-\sum_{l=1}^{\infty}\sum_{i=1}^{N}\prod_{j=1}^{D}\mathbb{P}\{x_{ij}(t_l+1) \in [lb_j, ub_j]\}\right)$$
$$\leq \exp\left(-\sum_{l=1}^{\infty}\sum_{i=1}^{N}\prod_{j=1}^{D}\varepsilon_{ij}\right) = 0.\ □$$

## Proof of lemma 6

Suppose that $h_{ij,\infty}(u_j) = 0$ holds. Just let $u_j > \frac{1}{3}\sum p_j$, since $h_{ij,t}(u_j)$ is a single-peak function symmetrical about $u_j = \frac{1}{3}\sum p_j$, $h_{ij,\infty}(u_j)$ becomes a single-peak function as well, and $\forall x \in (-\infty, \frac{2}{3}\sum p_j - u_j) \cup (u_j, +\infty), h_{ij,\infty}(x) = 0$, thus

$$\begin{aligned}
&\sigma^8 x_{ij}(t,\infty) \\
&= \int_{-\infty}^{+\infty} \left(u_j - \frac{1}{3}\sum p_j\right)^8 h_{ij,\infty}(u_j)\mathrm{d}u_j \\
&= \int_{\frac{2}{3}\sum p_j - u_j}^{u_j} \left(u - \frac{1}{3}\sum p_j\right)^8 h_{ij,\infty}(u_j)\mathrm{d}u_j \\
&\leq \left(u_j - \frac{1}{3}\sum p_j\right)^8 \int_{\frac{2}{3}\sum p_j - u_j}^{u_j} h_{ij,\infty}(u_j)\mathrm{d}u_j
\end{aligned}$$

$\leq \left(u_j - \frac{1}{3}\sum p_j\right)^8 < +\infty$, which contradicts with lemma 5! Thus, $h_{ij,\infty}(u_j) = 0$ is a false proposition.

## Proof of lemma 7

According to lemma 6, there exist a positive number $\varepsilon$ and an infinite subsequence $\{t_l\}, l = 1,2,\ldots$, such that $\forall l, h_{ij,t_l}(u_j) \geq \varepsilon$. Just suppose that $\left|ub_j - \frac{1}{3}\sum p_j\right| \geq \left|lb_j - \frac{1}{3}\sum p_j\right|$, thus $h_{ij,t_l}(ub_j) = \min_{x\in[lb_j.ub_j]} h_{ij,t_l}(x)$, and

$$\mathbb{P}\{x_{ij}(t_l+1) \in [lb_j, ub_j]\} = \int_{lb_j}^{ub_j} h_{ij,t_l}(u_j)\mathrm{d}u_j \geq (ub_j - lb_j)h_{ij,t_l}(ub_j) \geq (ub_j - lb_j)\varepsilon > 0$$

Accordingly, there exist a positive number $(ub_j - lb_j)\varepsilon$ and an infinite subsequence $\{t_l\}, l = 1,2,\ldots$, such that $\forall l, \mathbb{P}\{x_{ij}(t_l+1) \in [lb_j, ub_j]\} \geq (ub_j - lb_j)\varepsilon$, and $\lim_{T\to\infty} \lim_{t\to\lceil 0.05T\rceil} \mathbb{P}\{x_{ij}(t+1) \in [lb_j, ub_j]\} = 0$ is a false proposition. □

# Reference

[1] H. Zang, S. Zhang, K. Hapeshi, A review of nature-inspired algorithms, Journal of Bionic Engineering, 7 (2010) S232-S237, https://doi.org/10.1016/S1672-6529(09)60240-7.

[2] P. A. Vikhar, Evolutionary algorithms: A critical review and its future prospects, 2016 International Conference on Global Trends in Signal Processing, Information Computing and Communication (ICGTSPICC). IEEE, 2016, pp. 261-265, https://doi.org/10.1109/ICGTSPICC.2016.7955308.

[3] S. Mirjalili, S. M. Mirjalili, A. Lewis, Grey wolf optimizer, Advances in engineering software,

69 (2014) 46-61, https://doi.org/10.1016/j.advengsoft.2013.12.007.
[4] H. Wang, L. Shi, Y. Ni, Complementary and coordinated optimal dispatch solution for hybrid thermal-wind-PV power grid with energy storage, 8th Renewable Power Generation Conference (RPG 2019). IEEE, 2019, pp. 1-8, https://doi.org/10.1049/cp.2019.0544.
[5] F. Zhai, L. Shi, Solution proposal to the unit commitment problem incorporating manifold uncertainties, IET Generation, transmission & distribution, 14.21 (2020) 4763-4774, https://doi.org/10.1049/iet-gtd.2020.0805.
[6] S. Mirjalili, How effective is the Grey Wolf optimizer in training multi-layer perceptrons, Applied Intelligence, 43 (2015) 150-161, https://doi.org/10.1007/s10489-014-0645-7.
[7] M. R. Mosavi, M. Khishe, A. Ghamgosar. Classification of sonar data set using neural network trained by Gray wolf optimization, Neural Network World, 26.4 (2016) 393-415, https://doi.org/10.14311/NNW.2016.26.023.
[8] S. Zhang, Y. Zhou, Grey wolf optimizer based on Powell local optimization method for clustering analysis, Discrete Dynamics in Nature and Society, 2015 (2015) 1-17, https://doi.org/10.1155/2015/481360.
[9] V. Kumar, J. K. Chhabra, D. Kumar. Grey wolf algorithm-based clustering technique, Journal of Intelligent Systems, 26.1 (2017) 153-168, https://doi.org/10.1515/jisys-2014-0137.
[10] P. W. Tsai, T. T. Nguyen, T. K. Dao, Robot Path Planning Optimization Based on Multiobjective Grey Wolf Optimizer, in: J. S. Pan, J. W. Lin, C. H. Wang, X. Jiang, Genetic and Evolutionary Computing. ICGEC 2016. Advances in Intelligent Systems and Computing, vol 536. Springer, Cham, 2017, pp. 166-173, https://doi.org/10.1007/978-3-319-48490-7_20.
[11] C. Qu, W. Gai, M. Zhong, et al. A novel reinforcement learning based grey wolf optimizer algorithm for unmanned aerial vehicles (UAVs) path planning, Applied Soft Computing, 89 (2020) 106099. https://doi.org/10.1016/j.asoc.2020.106099.
[12] H. Faris, I. Aljarah, M. A. Al-Betar, et al., Grey wolf optimizer: a review of recent variants and applications, Neural Computing & Applications, 30 (2018) 413-435, https://doi.org/10.1007/s00521-017-3272-5.
[13] G. Negi, A. Kumar, S. Pant, et al., GWO: a review and applications, International Journal of System Assurance Engineering and Management, 12 (2021) 1-8, https://doi.org/10.1007/s13198-020-00995-8.
[14] N. M. Hatta, A. M. Zain, R. Sallehuddin, et al., Recent studies on optimisation method of Grey Wolf Optimiser (GWO): a review (2014–2017), Artificial Intelligence Review, 52 (2019) 2651-2683, https://doi.org/10.1007/s10462-018-9634-2.
[15] F. J. Solis, R. J. B. Wets, Minimization by Random Search Techniques, Mathematics of Operations Research, 3.1 (1981) 19-30, https://doi.org/10.1287/moor.6.1.19.
[16] J. Jgersküpper, Algorithmic analysis of a basic evolutionary algorithm for continuous optimization, Theoretical Computer Science, 379.3 (2007) 329-347, https://doi.org/10.1016/j.tcs.2007.02.042.
[17] G. Rudolph, Self-adaptive mutations may lead to premature convergence, IEEE Transactions on Evolutionary Computation, 5.4 (2001) 410-414, https://doi.org/10.1109/4235.942534.
[18] J. He, X. Yu, Conditions for the convergence of evolutionary algorithms, Journal of Systems Architecture, 47.7 (2001) 601-612, https://doi.org/10.1016/S1383-7621(01)00018-2.
[19] Y. Zhang, Some studies on the convergence and time complexity analysis of evolutionary algorithms, Diss. South China University of Technology, 2013.

[20] F. Van Den Bergh, A. P. Engelbrech, A Convergence Proof for the Particle Swarm Optimiser, Fundamenta Informaticae, 105.4 (2010), 341-374, https://doi.org/10.3233/FI-2010-370.
[21] F. Van Den Bergh, An analysis of particle swarm optimizers, Diss. University of Pretoria, 2007.
[22] Z. Hu, Q. Su, X. Yang, et al., Not guaranteeing convergence of differential evolution on a class of multimodal functions, Applied Soft Computing, 41 (2016) 479-487, https://doi.org/10.1016/j.asoc.2016.01.001.
[23] Z. Hu, S. Xiong, Q. Su, et al. Finite Markov chain analysis of classical differential evolution algorithm[J]. Journal of computational and applied mathematics, 268 (2014) 121-134, https://doi.org/10.1016/j.cam.2014.02.034.
[24] P. Feng, L. I. Xiao-Ting, Z. Qian, et al., Analysis of standard particle swarm optimization algorithm based on Markov chain, Acta Automatica Sinica, 39.4 (2013) 381-389, https://doi.org/10.1016/S1874-1029(13)60037-3.
[25] Z. Ren, J. Wang, Y. Gao, The global convergence analysis of particle swarm optimization algorithm based on Markov chain, Control Theory & Applications, 28.4 (2011) 462-466.
[26] G. Xu, G. Yu, Reprint of: On convergence analysis of particle swarm optimization algorithm, Journal of Computational and Applied Mathematics, 340 (2018) 709-717, https://doi.org/10.1016/j.cam.2018.04.036.
[27] D. Hu, X. Qiu, Y. Liu, et al., Probabilistic convergence analysis of the stochastic particle swarm optimization model without the stagnation assumption, Information Sciences, 547 (2021) 996-1007, https://doi.org/10.1016/j.ins.2020.08.072.
[28] J. Liu, X. Ma, X. Li, et al., Random convergence analysis of particle swarm optimization algorithm with time-varying attractor, Swarm and Evolutionary Computation, 61 (2021) 100819, https://doi.org/10.1016/j.swevo.2020.100819.
[29] Y. Tan, Y. Zhu, Fireworks Algorithm for Optimization, in: Y. Tan, Y. Shi, K. C. Tan, Advances in Swarm Intelligence. ICSI 2010. Lecture Notes in Computer Science, vol 6145. Springer, Berlin, Heidelberg, 2010, pp. 355-364, https://doi.org/10.1007/978-3-642-13495-1_44.
[30] J. Tong, The research and improvement of particle swarm optimization, Diss. Nanjing University of Information Science & Technology, 2008.
[31] Y. Ji, Two intelligent optimization algorithms their convergence analysis, Diss. Huazhong University of Science and Technology, 2011.
[32] J. Sun, B. Feng, W. Xu, Particle swarm optimization with particles having quantum behavior, in: Proceedings of the 2004 congress on evolutionary computation (IEEE Cat. No. 04TH8753), IEEE, 2004, pp. 325-331, https://doi.org/10.1109/CEC.2004.1330875.
[33] J. Sun, X. Wu, V. Palade, et al., Convergence analysis and improvements of quantum-behaved particle swarm optimization, Information Sciences, 193 (2012) 81-103, https://doi.org/10.1016/j.ins.2012.01.005.
[34] P. Chakraborty, S. Das, G. G. Roy, et al., On convergence of the multi-objective particle swarm optimizers[J]. Information Sciences, 181.8 (2011) 1411-1425, https://doi.org/10.1016/j.ins.2010.11.036.
[35] G. Xu, K. Luo, G. Jing, et al., On convergence analysis of multi-objective particle swarm optimization algorithm, European Journal of Operational Research, 286.1 (2020) 32-38, https://doi.org/10.1016/j.ejor.2020.03.035.
[36] Y. Zhang, D. Gong, X. Sun, et al. Adaptive bare-bones particle swarm optimization algorithm and its convergence analysis, Soft Computing, 18.7 (2014) 1337-1352,

https://doi.org/10.1007/s00500-013-1147-y.
[37] F. Pan, X. Hu, R. Eberhart, et al., An analysis of bare bones particle swarm, 2008 IEEE Swarm Intelligence Symposium, IEEE, 2008, pp: 1-5, https://doi.org/10.1109/SIS.2008.4668301.
[38] J. Liang, B. Qu, P. N. Suganthan, Problem definitions and evaluation criteria for the CEC 2014 special session and competition on single objective real-parameter numerical optimization, Computational Intelligence Laboratory, Zhengzhou University, Zhengzhou China and Technical Report, Nanyang Technological University, Singapore, 635 (2013) 490.
[39] X. S. Yang, Test problems in optimization, arXiv:1008.0549, 2010.